\newcommand{\eps}{\varepsilon}

\newcommand{\gam}{\gamma}
\newcommand{\fhi}{\varphi}
\newcommand{\ome}{\omega}

\newcommand{\sig}{\sigma}
\newcommand{\refa}[1]{{(\ref{#1})}}
\newcommand{\lam}{\lambda}
\newcommand{\Lam}{\Lambda}

\newcommand{\ov}{\overline}
\newcommand{\ii}{{\rm i}}

\newcommand{\Tf}{\widetilde \Phi_{\eps,N}}
\newcommand{\Beq}[2]{\begin{equation}\label{#1}{#2}\end{equation}}
\newcommand{\np}{\noindent}

\newcommand{\1} {{\frac{1}{ 2}}}

\newcommand{\R}{\mathbb{R}}
\newcommand{\Z}{\mathbb{Z}}
\newcommand{\C}{\mathbb{C}}
\newcommand{\N}{\mathbb{N}}
\newcommand{\Ca}{{\mathcal{C}}}
\newcommand{\La}{{\mathcal{L}}}
\newcommand{\G}{{\mathcal{G}}}
\newcommand{\Q} {\mathbb {Q}}
\newcommand{\T} {\mathbb {T}}
\newcommand{\Ba}{{\mathcal B}}

\newcommand{\B}{{\mathcal H}}
\newcommand{\I}{\int_{\T^2}}

\newcommand{\h}[1]{\vert #1 \vert_{H_1}}

\documentclass[12pt, a4 paper, twoside, openright]{amsart}
\usepackage{amsfonts}
\usepackage{latexsym}
 \usepackage{amsthm}
 \usepackage{amssymb}
\usepackage{amsopn}
\usepackage{layout}
\usepackage{fancyhdr}
\usepackage[english]{babel}
\usepackage{mathrsfs}
\usepackage[dvips]{graphicx}
\newcommand{\fig}[4]{\begin{figure}[!h]\begin{#1} 
\includegraphics[angle=0, width=#4\textwidth]{#2.eps} \caption{#3} 
\label{#2}\end{#1}
\end{figure}}
\theoremstyle{plain}

\newtheorem{pro}{Proposition}

\newtheorem{lem}{Lemma}[section]

\newtheorem{oss}{Remark}
\numberwithin{equation}{section}
\theoremstyle{definition}

\theoremstyle{remark}

\title[Quasi-periodic solutions of forced wave equations]{Quasi-periodic 
solutions of completely resonant forced 
wave equations}
\author{Massimiliano  Berti, Michela Procesi}

\thanks{ Massimiliano Berti, SISSA, via Beirut 2-4,
Trieste, Italy, \texttt{berti@sissa.it} \ . Michela Procesi,
Universit\`a di Roma $3$, Largo S. Leonardo Murialdo, Roma, Italy,
\texttt{procesi@mat.uniroma3.it \ .}}

\date{}
\begin{document}

\clearpage{\pagestyle{empty}\cleardoublepage}

\begin{abstract} 
We prove existence of quasi-periodic solutions with two 
frequencies of completely resonant, periodically 
{\it forced} nonlinear wave equations
with periodic spatial boundary conditions. We consider
both the cases the forcing frequency is: (Case A) a rational
number and (Case B) an irrational number.
\end{abstract}
\maketitle

{\small Keywords: Nonlinear Wave Equation, Quasi-Periodic Solutions,
Variational Methods, Lyapunov-Schmidt reduction,
Infinite Dimensional Hamiltonian Systems.\footnote{Supported by 
M.I.U.R. ``Variational Methods and Nonlinear
Differential Equations''.}
\\[1mm]
2000AMS subject classification: 35L05, 37K50, 58E05.}


\begin{section}{{\bf Introduction}}

We prove existence 
of small amplitude quasi-periodic solutions 
for completely resonant forced nonlinear wave equations
like
\Beq{e-0}{ 
\left\{\begin{aligned} & v_{tt}-v_{xx} +   f ( \ome_1 t , v ) = 0 \\ 
& v(t, x)= v(t, x + 2 \pi )\end{aligned}\right.} 
\np 
where the nonlinear forcing term 
$$
f ( \ome_1 t , v ) = a(\ome_1 t ) v^{2d-1} + O ( v^{2d} ), \qquad
\quad d > 1, \, 
d \in \N^+
$$
is $ 2\pi \slash \ome_1$-periodic in time. 
We shall consider both the cases
\begin{itemize}
\item
$ A) \ $ the forcing frequency $ \ome_1 \in \Q $ 
\item
$ B) \ $ the forcing frequency $ \ome_1 \in \R \setminus \Q
$.
\end{itemize}

Existence of 
periodic solutions for completely resonant forced wave equations
was first proved  in the pioneering papers \cite{R}, 
\cite{R1} (with Dirichlet boundary conditions)
if the forcing frequency is a rational number
($ \ome_1 = 1 $ in  \cite{R}-\cite{R1}).  
This requires to solve an infinite dimensional
bifurcation equation which lacks compactness property;
see \cite{bn}, \cite{Co}, \cite{BBi1}-\cite{BBi} and references therein 
for other results. 
If the forcing frequency is 
an irrational number 
existence of periodic solutions has been proved in 
\cite{PY}-\cite{Mc}: here the bifurcation equation is trivial but 
a ``small divisors problem'' appears. 
\\[1mm]
\indent
To prove existence of small amplitude
quasi-periodic
solutions for completely resonant PDE's like  (\ref{e-0}) one generally has to deal with a small divisor problem as well; however the main difficulty 
is to understand from which solutions of the linearized equation 
at $ v = 0 $,
$$
v_{tt} - v_{xx} = 0 \, ,
$$
quasi-periodic solutions branch-off: 
such linearized equation possesses only  $2\pi$-periodic solutions
$ q_+ (t+x) +$ $ q_-(t-x)$ where $ q_+( \cdot )$, $ q_-(\cdot ) $ are 
$ 2\pi$-periodic (completely resonant PDE). 

Here is the main difference w.r.t non-resonant PDE's for 
which a developed existence theory of 
periodic and quasi-periodic solutions has been established, 
see e.g. \cite{k1}, \cite{32}, \cite{cw}, \cite{p}, \cite{b3}
and references therein. 

For completely resonant autonomous PDE's,  
existence of periodic solutions has been proved in 
\cite{ls}, \cite{bp}, \cite{bb1}, \cite{bb2}, 
\cite{bb3}, \cite{gmp}, \cite{gp}, and 
quasi-periodic solutions with $2$-frequencies
have been recently obtained
in \cite{P1}-\cite{P2} for 
the specific nonlinearities $ f = u^3 + O(u^5) $.
Here the
bifurcation equation is solved by ODE methods.
\\[1mm]
\indent
In this paper we prove existence of quasi-periodic solutions 
with two frequencies $ \ome_1$, $\ome_2 $
for the completely resonant forced equation (\ref{e-0}) 
in both the two cases: Case $ A) $: $ \ome_1 \in \mathbb{Q} $;  
Case $ B) $: $ \ome_1 \in \R \setminus \Q $. 
\\[1mm]
\indent
The more interesting 
case  
is $ \ome_1 \in \Q $ (case A) when the 
forcing frequency $ \ome_1 $ enters in resonance 
with the linear frequency $1$. 
To find out from which solutions of the linearized equation 
quasi-periodic solutions of \refa{e-0} 
branch-off, requires to solve an infinite dimensional bifurcation 
equation which can not be solved 
in general by ODE
techniques
(it is a system of integro-differential equations). 
However, exploiting the variational nature of equation \refa{e-0} 
like in 
\cite{bb1}-\cite{bb2}, the bifurcation problem 
can be reduced to finding  critical points of
a suitable action functional which, in this case,  
possesses the infinite dimensional 
linking geometry \cite{BR}.    

\vskip8pt

\subsection{Main results} 
We  look for quasi-periodic solutions
$ v(t,x) $ of equation (\ref{e-0}) 
of the form
\Beq{e00}{ \left\{\begin{aligned}& v(t, x)=  u(\ome_1 t,\ome_2 t+x)\\
& u(\fhi_1+2k_1\pi,\fhi_2+2 k_2 \pi)=  u(\fhi_1,\fhi_2), \qquad 
\forall k_1, k_2 \in \Z
\end{aligned}\right.  }
with frequencies 
$$
\ome = (\ome_1, \ome_2) = (\ome_1, 1+ \eps) \, ,  
$$
imposing the frequency $ \ome_2 = 1 +\eps $ to be close 
to the linear frequency $ 1 $. 

Writing $ \partial_{tt} - \partial_{xx} = (\partial_t - \partial_x ) \circ
(\partial_t + \partial_x ) $ we get 
\begin{equation}\label{eqd}
\Big[ \ome_1 \partial_{\fhi_1} + (\ome_2-1) \partial_{\fhi_2} \Big]
\circ 
\Big[ \ome_1 \partial_{\fhi_1} + 
(\ome_2 + 1) \partial_{\fhi_2}  \Big] u + f( \varphi_1, u) = 0 
\end{equation}
and therefore 
\Beq{e01}{
\Big[
\ome_1^2 \partial^2_{\fhi_1}+(\ome_2^2-1)\partial^2_{\fhi_2} 
+ 2\ome_1\ome_2\partial_{\fhi_1}\partial_{\fhi_2} \Big] u(\fhi)
+  f(\fhi_1,u) = 0 \, .
  }

We assume that the forcing term $ f : \T \times \R \to \R $ 
$$
f (\fhi_1, u ) = a_{2d-1}(\fhi_1) 
u^{2d-1} + O( u^{2d} ), \qquad  d \in  \N^+ , \ d > 1 
$$  
is analytic in $ u $ but has only finite regularity
in $ \fhi_1 $. More precisely   

\medskip

\begin{itemize}
\item {\bf (H)} \ 
$ f ( \fhi_1, u ) := $ $ 
\sum_{k=2d-1}^\infty a_k(\fhi_1)u^k,$ $ d \in  \N^+ $, $ d > 1 $ 
and the coefficients $ a_k ( \fhi_1 ) \in H^1 ( \T ) $ verify,
for some $ r > 0 $, 
$\sum_{k=2d-1}^\infty | a_k |_{H^1} r^k  < \infty $.
The function $f(\fhi_1,u)$ is not identically constant in $\fhi_1$.
\end{itemize} 

\medskip

\np
We look for solutions $ u $ 
of (\ref{e01}) in the Banach space\footnote{Given  $ z \in \C $,
 $ z^* $ denotes its complex 
conjugate.} 
$$
\B_{\sig,s} :=  
\Big\{ 
u(\fhi) = 
\sum_{l \in \Z^2}  \hat u_l e^{\ii l \cdot \fhi}  \ \ : \ 
\hat u_l^* = \hat u_{-l} \ \  {\rm and} \ \    
| u |_{\sig, s} := \sum_{l \in \Z^2 } |\hat u_l | e^{|l_2 | \sig } 
[l_1]^{s}  < +\infty \Big\} 
$$
where $[l_1] := \max \{ |l_1| , 1 \} $ and $ \sigma > 0 $, $ s \geq 0 $.

The space $ \B_{\sig,s} $
is a Banach algebra 
with respect to multiplications of functions
(see Lemma \ref{banach} in the Appendix),
namely
$$
\quad u_1, u_2 \in \B_{\sig,s} \quad \Longrightarrow \quad 
u_1 u_2 \in \B_{\sig,s} \quad {\rm and} \quad 
|u_1 u_2|_{\sig,s} \leq C | u_1 |_{\sig,s} | u_2|_{\sig,s} \, .
$$

\np
We shall prove the following Theorems.

\medskip

\np
{\bf Theorem A.} {\it Let $\ome_1 = n\slash m \in \Q $. 
Assume that $ f $ satisfies assumption $(${\bf H}$)$ and
$ a_{2d-1}( \fhi_1 ) \neq 0 $, $\forall \fhi_1\in \T$.
Let $ \Ba_\gam $ be the uncountable\footnote{The proof 
that $ \Ba_\gam \cap (0, \eps_0 ) $ and $ \Ba_\gam \cap (- \eps_0,0 ) $
are both uncountable $ \forall \eps_0 > 0 $ is like in \cite{bp}.}
zero-measure Cantor set
$$  
\Ba_\gam 
:= \Big\{ \eps \in (-\eps_0, \eps_0)  \, : 
\quad |l_1 + \eps l_2 | > \frac{\gam}{|l_2| }\, , \;
\forall l_1, l_2\in \Z \setminus \{ 0 \} \Big\} 
$$ 
where $ 0 < \gam < 1 \slash 6 $.  

There exist constants
$ \overline \sig > 0 $, $ \overline s > 2 $, 
$ \overline{\eps} > 0 $,  $ \overline{C} > 0 $, 
such that  
$ \forall \eps \in \Ba_\gam $, 
$ |\eps | \gam^{-1} \leq \bar{\eps}\slash m^2$,  
there exists a classical solution $ u ( \eps, \fhi ) \in 
\B_{\overline{\sig}, \overline{s}} $ 
of (\ref{e01}) with 
$(\ome_1, \ome_2 ) = (n\slash m, 1+ \eps ) $
satisfying 
\begin{equation}\label{ba12}
\Big| u(\eps,\fhi) - |\eps |^{ \frac{1}{2(d-1)}} \bar q_\eps (\fhi)
\Big|_{\overline{\sig},\overline{s}} \leq 
\overline{C}  \frac{m^2|\eps |}{\gam \, \ome_1^3} 
|\eps |^{ \frac{1}{2(d-1)}} 
\end{equation}
for an appropriate function 
$ \bar q_\eps \in \B_{\sig,s} \setminus \{ 0 \}$ of the form
$ \bar q_\eps (\fhi)= $ $\bar q_+(\fhi_2)+ $ 
$\bar q_-(2m\fhi_1-n\fhi_2)$.

As a consequence, equation (\ref{e-0})
admits 
the quasi-periodic solution  
$ v(\eps, t, x) := $ 
$ u(\eps, \ome_1 t, x+ \ome_2 t ) $
with two frequencies
$ (\ome_1, \ome_2 ) = $ $ (n\slash m, 1 + \eps )$ and 
the map $ t \to v(\eps, t, \cdot ) \in 
H^{\overline{\sig}} ( \T )$
has the form\footnote{We denote  $ H^{\sig } (\T) := \{  
u( \fhi ) = \sum_{l \in {\Z}} \hat u_l e^{\ii l \fhi}  \ : \
\hat u_l^* = \hat u_{-l} \, , \   
|u|_{H^{\sig}(\T)} := 
\sum_{l\in \Z} |\hat u_l |  e^{\sig |l|} < + \infty \} $.} 
$$ 
\Big| v(\eps, t, x ) -
|\eps |^{\frac{1}{2(d-1)}} \Big[ \bar q_+ ( x + ( 1 + \eps ) t ) +
\bar q_- ( (1 - \eps) nt - nx ) \Big] \Big|_{H^{\overline{\sig}}( \T )} =  
O \Big( \frac{m^2}{\gam \, \ome_1^3} 
|\eps |^{ \frac{2d-1}{2(d-1)}} \Big) \, .
$$}

\medskip

\np
At the first order the quasi-periodic solution $v (\eps, t , x) $ 
of equation \refa{e-0}
 is the superposition of two waves 
traveling
in opposite directions 
(in general, both components $q_+ $, $q_-$ 
are non trivial). 

The bifurcation of quasi-periodic solutions 
looks quite different if $\ome_1 $ is irrational.

\medskip

\np
{\bf Theorem B.} {\it \label{TB} 
Let $ \ome_1 \in \R \setminus \Q $. 
Assume that $ f $ satisfies assumption $(${\bf H}$)$, 
$ \int_0^{2\pi} a_{2d-1} (\fhi_1) d \fhi_1 $ $ \neq 0 $ and 
$ f (\fhi_1,u)\in H^s(\T) $, $ s \geq 1 $, for all $ u $. 

Let $ \Ca_\gam \subset D \equiv (-\eps_0,\eps_0) 
\times(1,2)$ be the uncountable zero-measure Cantor set\footnote{See 
Lemma \ref{Cgam}.}
\begin{equation}\label{cgam}
\Ca_\gam := \left\{\begin{aligned}(\eps,\ome_1)\in D:\;\; 
& \ome_1\notin \Q\,, \quad \frac{\ome_1}{\ome_2}\notin \Q \,,\quad  
|\ome_1 l_1 +\eps l_2|>\frac{\gam}{|l_1| + |l_2|}\,,\\ 
& |\ome_1 l_1 + (2 + \eps  ) l_2| >
\frac{\gam}{|l_1| + |l_2|}\,, \ \forall \; l_1,l_2 \in \Z\setminus 
\{0\}\quad\end{aligned}\right\}.
\end{equation}
Fix any $ 0 < \overline{s} < s - 1 \slash 2 $. 
There exist positive constants $ \overline{\eps} $, 
$ \overline{C } $, $ \overline{\sig} > 0 $, 
such that, 
$ \forall ( \eps, \ome_1 ) \in \Ca_\gam $ with
$ |\eps | \gamma^{-1} < \overline{\eps} $ and 
$ \, \eps \int_0^{2\pi} a_{2d-1} (\fhi_1) \, d \fhi_1> 0$, there exists 
a nontrivial solution 
$ u (\eps, \fhi ) \in \B_{\overline{\sig},\overline{s}} $
of equation (\ref{e01}) with $(\ome_1, \ome_2 ) = (\ome_1, 1+\eps)$
satisfying
\begin{equation} \label{sB}
\Big|u(\eps, \fhi) -  |\eps|^{\frac{1}{2(d-1)}}
{\bar q}_\eps ( \fhi_2) \Big|_{\overline{\sig}, \overline{s}} 
\leq  \ \overline{C } 
\frac{| \eps |}{\gamma}  |\eps|^{\frac{1}{2(d-1)}}
\end{equation}
for some function $ {\bar q}_\eps (\fhi_2) \in 
H^{\overline{\sig}} (\T) \setminus \{ 0 \} $. 

As a consequence, 
equation (\ref{e-0}) admits the non-trivial quasi-periodic solution  
$ v(\eps, t, x) := $ $ u(\eps, \ome_1 t, x + \ome_2 t  )$ with
two frequencies 
$ (\ome_1, \ome_2 ) = $ $ (\ome_1, 1 + \eps )$ and the 
map $ t \to v(\eps, t , \cdot ) \in H^{\overline{\sig}}( \T ) $ 
has the form
$$ 
\Big| v(\eps, t, x ) - |\eps|^{\frac{1}{2(d-1)}}
{\bar q}_\eps ( x + ( 1 + \eps ) t )  
\Big|_{H^{\overline{\sig}}( \T )}
= O \Big( \gamma^{-1} |\eps |^{\frac{2d-1}{2(d-1)}} \Big)  .
$$}

\begin{oss} 
Imposing in the definition of $ \Ca_\gamma $ 
the condition $ \ome_1 \slash \ome_2 = \ome_1 \slash (1+\eps ) \in 
\Q $ 
we obtain, by Theorem B the
existence of periodic solutions of equation  (\ref{e-0}).
They are reminiscent, in this 
completely resonant context, of the Birkhoff-Lewis periodic orbits
with large minimal period accumulating at the origin, 
see \cite{BaB}, \cite{BBV}. 
\end{oss}

\begin{oss} {\bf (Non existence)}
In Theorem  B, existence of quasi-periodic solutions 
could follow by other hypotheses on $ f $, see remark 
\ref{othersol}.
However the hypothesis that the leading term in the nonlinearity $ f $ 
is an odd power of $ u $ is not of purely technical nature.
If $ f(\fhi_1, u ) = a ( \fhi_1) u^D $ with 
$ D $ even and $ \int_0^{2\pi} a ( \fhi_1) \, d \fhi_1 \neq 0 $, then,
$ \forall R > 0 $ 
there exists  $ \eps_0 > 0 $ such that 
$ \forall \sig \geq 0$, $ \overline{s} > s - 1\slash 2 $, 
$ \forall (\eps, \ome_1) \in {\mathcal C}_\gamma$
with $|\eps | < \eps_0 $, equation (\ref{e01}) does not possess  
solutions $ u \in \B_{\sig,\overline{s}}$ in the ball 
$ |u|_{\sig, \overline{s}}\le 
R |\eps|^{1\slash (D-1)}$, see Proposition \ref{nonex}.
\end{oss}

To prove Theorems A-B, 
instead of looking for solutions of equation \refa{e01} in a shrinking 
neighborhood of $0$, it is a convenient devise to perform the rescaling
$$
u \to \delta  u
\qquad \ {\rm with } \ \qquad 
\delta := |\eps |^{1 \slash 2(d-1)}  
$$ 
enhancing the relation between the amplitude $ \delta $ and the 
frequency $ \ome_2 = 1 + \eps $. We obtain the equation 
\begin{equation}\label{e1}
{\mathcal L}_\eps u + \eps f(\fhi_1,u, \delta ) = 0  
\end{equation}
where, see \refa{eqd},  
\begin{eqnarray*}
{\mathcal L}_\eps & := & 
\Big[ \ome_1 \partial_{\fhi_1} + \eps \partial_{\fhi_2} \Big]
\circ 
\Big[ \ome_1 \partial_{\fhi_1} + (2+ \eps )  \partial_{\fhi_2} \Big]   
\\
& = & 
\Big[ \ome_1^2 \partial^2_{\fhi_1} + 
2 \ome_1 \partial_{\fhi_1}\partial_{\fhi_2} \Big] +
\eps \Big[ (2 + \eps ) \partial^2_{\fhi_2} 
+ 2 \ome_1 \partial_{\fhi_1}\partial_{\fhi_2} \Big] 
\end{eqnarray*}
and 
\begin{equation}\label{fd}
f(\fhi_1,u , \delta ) := {\rm sign}(\eps) 
\frac{f(\fhi_1, \delta u )}{\delta^{2(d-1)}} =
{\rm sign}(\eps) \Big( 
a_{2d-1} ( \fhi_1 )u^{2d-1} + \delta a_{2d} (\fhi_1 ) u^{2d} + \ldots \Big)
\end{equation}
and sign$(\eps) := 1 $ if $\ome_2 > 1$ and 
sign$(\eps) := - 1 $ if $ \ome_2 < 1 $.

\vskip5pt

To find solutions of equation (\ref{e1}) 
we shall apply the Lyapunov-Schmidt decomposition method which leads to 
solve separately a ``range equation'' and a ``bifurcation equation''.

In order to solve the range equation (avoiding small divisor problems)
we restrict $ \eps $ to the uncountable zero-measure 
set $ {\mathcal B}_\gamma $ for Theorem A, resp. 
$ (\eps, \ome_1) \in {\mathcal C}_\gamma $ for Theorem B, 
and we apply the Contraction Mapping Theorem; similar
non-resonance conditions have been employed e.g. in 
\cite{ls}, \cite{bp},
\cite{bb1}-\cite{bb2}, \cite{Mc}, \cite{P1}.

To solve the {\it infinite} dimensional 
bifurcation equation we proceed in different ways 
in  case A) and case B). 

As already said, 
in case A) we follow the variational approach of \cite{bb1},\cite{bb2} 
noting that the bifurcation equation is the Euler-Lagrange equation
of a ``reduced action functional'' which 
turns out to have the geometry of the infinite dimensional linking theorem
of Benci-Rabinowitz \cite{BR}. 
However we can not directly apply the linking theorem because 
the reduced action functional is defined only in a ball centered at the origin 
(where 
the range equation is solved). 
Moreover the infinite dimensional linking theorem of \cite{BR} 
requires the compactness of the gradient of the functional, property
which is not preserved by extending the functional in the whole 
infinite dimensional space.

In order to overcome these difficulties we perform
a further finite dimensional reduction 
of Galerkin type inspired to \cite{bb3} 
on a subspace of dimension $N$, with $N$ 
large but independent of $ \eps $, see 
the equations \refa{q1}-\refa{q2}-\refa{p}. 

We shall have to solve the \refa{q2}-\refa{p} equations in a 
sufficiently large domain of $ q_1 $ (Lemma \ref{L2}), 
consistent with the $| \cdot |_{H^1}$ 
bounds on the solution $ q_1 $ of the bifurcation equation 
that can be obtained by the variational
arguments, see Lemma \ref{L28}. 

Another advantage of this method is that allows 
to prove the analiticity of the solution $ u $  
in the variable $ \fhi_2 $.

In case B) the 
bifurcation equation could be solved through 
variational methods as in case A). However there is a simpler 
technique available. The bifurcation equation reduces,
in the limit $ \eps\rightarrow 0$, 
to a super-quadratic Hamiltonian system 
with one degree of freedom. We  prove existence of  
a non-degenerate solution by
phase-space analysis. Therefore it can be continued 
by the Implicit function Theorem 
to a solution of the complete bifurcation equation 
for $ \eps $ small.  

\end{section}

\bigskip

The paper is organized as follows. For simplicity of exposition 
we prove first Theorem A in the case $ \ome_1 = 1 $. 
We deal with the general case $ \ome_1 = \frac{n}{m} \in \Q $
at the end of section \ref{caseA}. In section \ref{caseB}
we prove Theorem B.

\bigskip

{\bf Acknowledgments:} The authors thank Luca Biasco and Philippe Bolle
for useful comments.
Part of this paper was written when the second author was at SISSA.

\begin{section}{{\bf Case A: $ \ome_1 \in \Q $}}\label{caseA}

Equation (\ref{e1}) becomes, for $ \ome_1 = 1 $ 
\begin{equation}\label{l1}
{\mathcal L}_\eps u + \eps f(\fhi_1,u, \delta ) = 0  
\end{equation}
where 
\begin{eqnarray*}
{\mathcal L}_\eps & := & 
\Big[ \partial_{\fhi_1} + \eps \partial_{\fhi_2} \Big]
\circ 
\Big[ \partial_{\fhi_1} + (2+\eps) \partial_{\fhi_2} \Big]   \\
& = & 
\Big[ \partial^2_{\fhi_1} + 2 \partial_{\fhi_1}\partial_{\fhi_2} \Big] +
\eps \Big[ (2 + \eps ) \partial^2_{\fhi_2} 
+ 2 \partial_{\fhi_1}\partial_{\fhi_2} \Big] 
\equiv L_0 + \eps L_1 \, .
\end{eqnarray*}

To fix notations we shall prove 
Theorem A in the case $a_{2d-1}(\fhi_1)>0 $ and $\eps > 0 $, i.e. 
sign$(\eps) > 0 $.

By the assumption {\bf (H)} on the nonlinearity $ f $ 
and by the Banach algebra property of $ \B_{\sig,s} $ 
the Nemitskii operator
$$  
u \to f(\fhi_1, u , \delta ) 
\in \  C^\infty  ( B_\rho, \B_{\sig,s} )\, ,  \qquad \qquad
0 < s < \frac12   
$$ 
where $ B_\rho $ is the ball of radius 
$ \rho \delta^{-1} $ in  $\B_{\sig,s} $ and $\rho$ is connected 
to the analiticity radius $ r $ of $ f $
(note that since $ a_k ( \fhi_1 ) \in H^1 (\T) $ then $ a_k ( \cdot ) 
\in {\mathcal H}_{\sig,s}$, $ \forall \sig > 0 $,
$ 0 < s <  1 \slash 2 $).

\medskip

Equation (\ref{l1}) is the Euler-Lagrange equation of the 
{\em Lagrangian action functional} 
$ \Psi_\eps 
\in C^1(\B_{\sig,s},\R)$ 
defined by
\begin{eqnarray*} 
\Psi_\eps (u) & := & \I \1 ( \partial_{\fhi_1} u )^2 + 
( \partial_{\fhi_1} u) ( \partial_{\fhi_2} u ) + 
\frac{  \eps (2 + \eps) }{2} ( \partial_{\fhi_2} u )^2
+  \eps( \partial_{\fhi_1} u) ( \partial_{\fhi_2} u )  
-  \eps F( \fhi_1 , u, \delta ) \\
& \equiv &  \Psi_0 (u) +  \eps \Gamma (u, \delta ) 
\end{eqnarray*}
where $ F ( \fhi_1 , u,\delta  ) := \int_0^u 
f ( \fhi_1 , \xi , \delta )d\xi \ $  and
\begin{eqnarray*}
\Psi_0 (u) &:= & \I \1 ( \partial_{\fhi_1}u)^2 
+ (\partial_{\fhi_1} u ) ( \partial_{\fhi_2}u)  \label{psi0} \\
\Gamma (u, \delta ) & := & 
\I \frac{(2 + \eps )}{2} ( \partial_{\fhi_2} u )^2
+ ( \partial_{\fhi_1} u) ( \partial_{\fhi_2} u ) 
- F( \fhi_1 , u, \delta ) \, .   
\label{gamma}
\end{eqnarray*}
\\[1mm]
To find critical points of $ \Psi_\eps $ 
we perform a variational Lyapunov-Schmidt reduction
inspired to \cite{bb1}-\cite{bb2}, see also \cite{AB}.

\smallskip

\subsection{The Variational Lyapunov-Schmidt Reduction}

The unperturbed functional $ \Psi_0 : \B_{\sig,s} \to \R $ 
possesses 
an infinite dimensional linear space $ Q $ of critical points 
which are the solutions $ q $ 
of the equation 
$$ 
L_0 q = 
\partial_{\fhi_1} \Big( \partial_{\fhi_1} + 2 \partial_{\fhi_2} \Big) q 
= 0 \ .
$$ 
The space $ Q $  can be written as
$$
Q  =   \Big\{ q = \sum_{l \in \Z^2 }
\hat q_l e^{\ii l \cdot \fhi } \in \B_{\sig,s} \  | \ \ 
\hat q_l = 0 \ \ \ {\rm for } \ \ \ l_1 (l_1 + 2 l_2) \neq 0 \Big\} \ .
$$
In view of the variational argument 
that we shall use to solve the bifurcation equation
we split $ Q $ as 
$$
Q = Q_+\oplus Q_0 \oplus Q_-
$$
where\footnote{$ H^{\sig}_0 ( \T ) $ denotes the functions  
of $ H^{\sig} (\T)$ with zero average.  
$ H^{\sig,s} (\T) :=$ $\{  
u( \fhi ) = \sum_{l \in {\Z}} \hat u_l e^{\ii l \fhi} \ : $
$ \hat u_l^* = \hat u_{-l}, $ 
$ |u|_{H^{\sig,s}(\T)} :=$ 
$\sum_{l\in \Z} |\hat u_l |  e^{\sig |l|}[l]^s < + \infty \} $
and $ H^{\sig,s}_0 (\T) $ its functions with zero average.} 
\begin{eqnarray*}
Q_+ & := & \Big\{ q \in Q \ : \ \hat q_l = 0 \ {\rm for} \ 
l \notin \Lam_+ 
\Big\} = \Big\{ q_+ := q_+ (\fhi_2 ) \in H^{\sig}_0 (\T) \Big\} \\
Q_0 &:= & \Big\{  q_0 \in \R \Big\} \\
Q_- &: = & \Big\{ q \in Q \ : \ \hat q_l = 0 \ {\rm for} \
l \notin \Lam_- 
\Big\} = \Big\{ q_- := 
q_- (2 \fhi_1 -\fhi_2 ), \ q_- (\cdot ) \in H^{\sig,s}_0 (\T)\Big\}
\end{eqnarray*}
and 
\begin{equation}\label{L+-}
\Lam_+ := \Big\{ l \in \Z^2 \ : \ l_1 = 0, \ l \neq 0 \Big\}, 
\quad 
\Lam_- := \Big\{ l \in \Z^2 \ : \ l_1 + 2 l_2 = 0, \ l \neq 0 \Big\} \, .
\end{equation}
We shall also use in $ Q $ the norm 
$$
| q |_{H^1}^2 = |q_+|_{H^1(\T)}^2 + q_0^2 + |q_-|_{H^1(\T)}^2 \,
\sim \,
\sum_{l\in \Lambda_- \cup \{0\} \cup \Lambda_+ } \hat q_l^2 (|l|^2 +1) \, .
$$
We decompose the space 
$ \B_{\sig,s} = Q \oplus P $
where 
\begin{eqnarray*}
P & := & 
\Big\{ p = \sum_{l \in \mathbb{Z}^2}
{\hat p}_l e^{\ii l \cdot  \fhi } \in \B_{\sig,s}  \ | \ \  
{\hat p}_l = 0 \ \ \ {\rm for} \  \ \ l_1 (2 l_2 + l_1) = 0 \Big \}  \ .
\end{eqnarray*}
Projecting equation  
(\ref{l1}) onto the closed subspaces $ Q $ and $ P $, 
setting $ u = q +  p \in \B_{\sig,s} $ with $ q \in Q $ and 
$ p \in P $, we obtain
$$
\left\{\begin{aligned} & L_1 [q] + \ \Pi_Q f(\fhi_1, q + p, \delta )  = 0  
\qquad \quad (Q)\\ 
&
{\mathcal L}_\eps[p] + \eps \Pi_{P}  f(\fhi_1, q+p , \delta ) = 0  
\qquad \quad (P)
\end{aligned} \right.
$$
where $ \Pi_Q : \B_{\sig,s} \to Q $, $ \Pi_P : \B_{\sig,s} \to P $
are the projectors respectively onto $ Q $ and $ P $.
\\[1mm]
\indent
In order to prove analiticity of the solutions and
to highlight the compactness of the problem we perform
a {\it finite} dimensional Lyapunov-Schmidt  reduction, introducing 
the decomposition 
$$
Q = Q_1 \oplus Q_2
$$
where 
\begin{eqnarray*}
Q_1 := Q_1(N) := \Big\{ \ q = 
\sum_{ |l| \leq N }  \hat q_l e^{\ii l \cdot \fhi } \in Q \Big\} \, ,
\ Q_2 := Q_2(N) := \Big\{ q = 
\sum_{ |l| > N }  \hat q_l e^{\ii l \cdot \fhi } \in Q \Big\} \ .
\end{eqnarray*} 
Setting $q = q_1 + q_2 $ with $ q_1 \in Q_1 $ and $ q_2 \in Q_2 $, we 
finally get 
\begin{eqnarray}
\qquad L_1 [q_1] + \Pi_{Q_1} 
\Big[ f(\fhi_1, q_1 + q_2 + p, \delta ) \Big] = 0  
& \iff & d \Psi_\eps ( u )[h] = 0 \ \  \forall h \in Q_1  
\quad (Q_1)  \label{q1}
\\
\qquad L_1 [q_2] 
+ \Pi_{Q_2} \Big[ f(\fhi_1, q_1 + q_2 + p, \delta ) \Big] = 0  
& \iff & d \Psi_\eps ( u )[h] = 0 \ \  \forall h \in Q_2  
\quad (Q_2) \label{q2}
\\
\qquad {\mathcal L}_\eps[p] 
+ \eps \Pi_{P} \Big[ f(\fhi_1, q_1 + q_2 + p, \delta )\Big] = 0  
& \iff &  d \Psi_\eps ( u )[h] = 0 \ \  \forall h \in P \quad \ (P) \label{p}
\end{eqnarray}
where $ \Pi_{Q_i} : \B_{\sig,s} \to Q_i $ 
are the projectors onto $ Q_i $ ($ i = 1, 2 $).
\\[1mm]
\indent
We shall solve first the $(Q_2)$-$(P)$-equations 
for all $ | q_1 |_{H^1} \leq 2R $, 
provided $ \eps $ belongs to a suitable Cantor-like set,
$ |\eps | \leq \eps_0 (R) $ is sufficiently small and 
$ N \geq  N_0 (R) $ is large enough (see Lemma \ref{L2}).

Next we shall solve the $(Q_1)$-equation by means of 
a variational linking argument, 
see subsection \ref{lin}.

\medskip

\begin{subsection}{The $(Q_2)$-$(P)$-equations}
We first prove that
$ {\mathcal L}_\eps $ restricted to $ P $ has 
a bounded inverse 
when $ \eps $ belongs to the uncountable
zero measure set 
$$  
\Ba_\gam 
:= \Big\{ \eps \in (- \eps_0, \eps_0)  \, : 
\quad |l_1 + \eps l_2 | > \frac{\gam}{|l_2| }\, , \;
\forall l_1, l_2\in \Z \setminus \{ 0 \} \Big\}
$$
where $ 0 < \gam < 1 \slash  6 $.  
$ \Ba_\gam $ accumulates at $0$ both from the right and from the left, see
\cite{bp}.
\\[1mm]
\indent
The operator $ {\mathcal L}_\eps $ is diagonal in the 
Fourier basis 
$ \{ e^{\ii l \cdot \fhi } $ , 
$ l \in \mathbb{Z}^2 \} $ with  eigenvalues
$ D_l :=$ $( l_1 +\eps l_2)( l_1 + (2  + \eps  ) l_2)$.

\begin{lem}\label{L1} For $ \eps \in \Ba_\gam $ 
the eigenvalues $ D_l $ of $ {\mathcal L}_\eps $ 
restricted to $P$, satisfy
$$
|D_l| = \Big| l_1 +\eps l_2  \Big| \, \Big| 
(l_1 + 2 l_2) + \eps l_2 \Big|
> \gam 
\qquad \forall l_1 \neq 0, \ l_1 + 2l_2 \neq 0  \ .
$$
As a consequence 
the operator $ {\mathcal L}_\eps : P \to P $ 
has a bounded inverse ${\mathcal L}_\eps^{-1}$ 
satisfying
\Beq{b1}{ \quad \Big| {\mathcal L}_\eps^{-1}[h] \Big|_{\sig,s} \leq 
\frac{|h|_{\sig,s}}{\gamma}, \qquad \forall h \in P \, .}
\end{lem} 

\begin{proof} Denoting by $[x]$ the nearest integer  
close to $ x $ and $\{x\}=x-[x]$, we have that 
$ D_l > 1 $ if both $ l_1 \neq- [\eps l_2]$ and 
$ l_1 + 2 l_2\neq -[\eps l_2] $. 
If $ l_1 = - [\eps l_2] $ then 
$$
|D_l|\geq \frac{\gam}{|l_2|} (|2 l_2|-\{\eps l_2\})\geq \gam. 
$$
In the same way if  $ l_1+2 l_2= -[\eps l_2] $ we have 
$ |D_l|\geq$ $\frac{\gam}{|l_2|} (|2l_2|-\{\eps l_2\})\geq$ $ \gam. $
\end{proof}

\begin{lem}\label{L1inv} 
The operator $ L_1 : Q_2 \to Q_2 $ 
has bounded inverse $ L_1^{-1} $ which satisfies
\Beq{l1invN}{\Big| L_1^{-1}[h]  \Big|_{\sig,s} \leq \frac{|h|_{\sig,s}}{N^2}. }
\end{lem} 

\begin{proof} 
$L_1$ is diagonal in the Fourier basis of $Q$:
$ e^{\ii l \cdot \fhi} $ with
$l \in \Lambda_+ \cup \{0\} \cup \Lambda_- $ (recall \refa{L+-})
with eigenvalues
\begin{equation}\label{L1-1}
d_l= (2+\eps)l_2^2 \ \ {\rm if} \ \ 
l_1=0 \quad {\rm and} 
\quad d_l = (-2+\eps)l_2^2 \ \ {\rm if} \ \  l_1 + 2l_2=0 \, .
\end{equation}
The eigenvalues of $ L_1 $ restricted to $ Q_2 ( N ) $ 
verify  $ |d_l| \geq (2 - \eps) N^2 $ and
\refa{l1invN} holds.
\end{proof}


\medskip

Fixed points of the nonlinear operator ${\mathcal G} : 
Q_2 \oplus P \to Q_2 \oplus P $
defined by 
$$
{\mathcal G} (q_2, p ; q_1) 
:= \Big(- L_1^{-1} \Pi_{Q_2}  f(\fhi_1, q_1 + q_2 + p, \delta  ), \  
- \eps {\mathcal L}_\eps^{-1} \Pi_P f(\fhi_1 , q_1 + q_2 + p,\delta )  \Big)
$$
are solutions of the $(Q_2)$-$(P)$-equations.

Using the Contraction Mapping Theorem we can prove:

\begin{lem}\label{L2}
{\bf (Solution of the ($Q_2$)-($P$) equations)}
$ \forall R  > 0 $ there exist an integer $  N_0 (R) 
\in {\mathbb N}^+ $ and positive constants $ \eps_0 (R)> 0 $, 
$ C_0 ( R ) > 0 $ such that:
\begin{equation} \label{hp}
\forall  |q_1|_{H^1} \leq 2 R \,,\;\;
\forall \eps \in B_\gam, \ | \eps | 
\gamma^{-1} \leq \eps_0 (R) \,,\;\; \forall  
 N\ge N_0 ( R )  \ : \ 
0 \leq \sig N \leq 1\,, 
\end{equation}

\np
there exists a unique solution $ (q_2 (q_1), p (q_1)) := $ 
$( q_2(\eps,N,  q_1 ), p(\eps, N,  q_1 ) )
\in $ $ Q_2 \oplus P $ 
of the $(Q_2)$-$(P)$ equations
satisfying
\Beq{co1}{
|q_2 (\eps,N, q_1)|_{\sig,s} \leq \frac{C_0(R)}{N^2}, \quad
| p (\eps ,N, q_1 ) |_{\sig,s} \leq C_0(R) | \eps | \gamma^{-1} \, .
}
Moreover the map $ q_1 \rightarrow (q_2 (q_1), 
p (q_1)) $ is in $ C^1 ( B_{2R}, Q_2 \oplus P)$ and 
\Beq{co2}{
\Big| p'(q_1)[h] \Big|_{\sig,s} \leq C_0(R) |\eps|\gamma^{-1} |h|_{H^1} \, , 
\quad 
\Big| q_2'(q_1)[h] \Big|_{\sig,s} \leq \frac{C_0(R)}{N^2} |h|_{H^1} \quad
\forall h \in Q_1 \, .}
\end{lem}

\begin{proof} In the Appendix.
\end{proof}

\end{subsection}

\medskip

\begin{subsection}{The $(Q_1)$-equation}\label{sub:Q1}

Once the $(Q_2)$-$(P)$-equations have been solved by 
$(q_2 (q_1),$ $ p(q_1)) \in $ $ Q_2 \oplus P $
there remains the finite dimensional $(Q_1)$-equation 
\begin{equation}\label{PQ2p}
L_1 [q_1] + \Pi_{Q_1}  
f(\fhi_1, q_1 + q_2(q_1) + p (q_1),\delta ) 
= 0 \ .   
\end{equation}

\medskip

The geometric interpretation of the construction of $ (q_2(q_1)$, 
$ p (q_1)) $ is that on the finite 
dimensional sub-manifold $ Z \equiv $
$ \{ q_1 + q_2 (q_1) + p (q_1) \ : \ |q_1| < 2R \}$,  
diffeomorphic to the ball
$$
B_{2R}:= \{ q_1 \in Q_1 \ : \ 
|q_1|_{H^1} < 2 R \}, 
$$ 
the partial derivatives of the action
functional $ \Psi_\eps $ with respect to the variables $ (q_2, p) $
vanish.  
We claim that at a critical point of $ \Psi_\eps $ 
restricted to $ Z $, also the partial derivative
of $ \Psi_\eps $ w.r.t. the variable $ q_1 $ vanishes and therefore 
that such point is critical also for the non-restricted functional
$ \Psi_\eps : \B_{\sig,s} \to \R  $.

Actually the bifurcation equation (\ref{PQ2p})
is the Euler-Lagrange equation 
of the reduced Lagrangian action functional 
$$
\Phi_{\eps,N} : B_{2R} \subset Q_1 
\to \mathbb{R},  \qquad \quad
\Phi_{\eps,N} (q_1 ) := \Psi_\eps (q_1 + q_2(q_1) + p(q_1) ) \ .
$$

\begin{lem}\label{L24}
$ \Phi_{\eps,N} \in $ $ \mathcal{C}^1( B_{2R}, \mathbb{R})$ and a
critical point $ q_1 \in B_{2R} $ 
of $ \Phi_{\eps,N} $ is a solution of the bifurcation equation (\ref{PQ2p}).
Moreover $ \Phi_{\eps,N} $ can be written as
\begin{equation}\label{r1}
\Phi_{\eps,N} ( q_1 ) =  {\rm const} + \eps \Big( 
\Gamma ( q_1 )  + {\mathcal R}_{\eps,N}(q_1) \Big)
\end{equation}
where
$$ 
\Gamma(q_1) :=  \I \frac{(2 + \eps )}{2} ( \partial_{\fhi_2} q_1 )^2
+ ( \partial_{\fhi_1} q_1) ( \partial_{\fhi_2} q_1) 
- a_{2d-1}(\fhi_1) \frac{q_1^{2d}}{2d}  
$$  
$$
\begin{aligned}{\mathcal R}_{\eps,N}(q_1) := &
\I F(\fhi_1, q_1 , \delta = 0 ) - F(\fhi_1, q_1+ q_2 (q_1) + p (q_1),\delta) 
\\
& + \frac{1}{2} f(\fhi_1, q_1+ 
q_2(q_1) + p(q_1),\delta) (q_2 (q_1) + p (q_1))  \end{aligned}
$$
and, for some positive constant $C_2 (R) \geq C_1 (R) $,
\begin{eqnarray}\label{Restun}
| {\mathcal R}_{\eps,N} (q_1) | & \leq & 
C_2(R) \Big(\delta+ |\eps|\gam^{-1} + \frac{1}{N^2} \Big) \\
\label{Rest1un}
\Big| {\mathcal R}_{\eps,N}' (q_1)[h] \Big|  & \leq &
C_2(R) \Big(\delta+ |\eps|\gam^{-1} + \frac{1}{N^2} \Big) | h |_{H^1}, 
\qquad \forall h \in Q_1 \, . 
\end{eqnarray}
\end{lem}

\begin{proof}
In the Appendix. 
\end{proof}

The problem of finding non-trivial solutions
of the $Q_1$-equation is reduced to finding 
non-trivial critical points of the reduced action functional
$  \Phi_{\eps,N} $ in $ B_{2R} $.

By (\ref{r1}), this is equivalent to find critical
points of the rescaled functional 
(still denoted $ \Phi_{\eps,N} $ and called the reduced
action functional)
\begin{equation}\label{redf} 
\Phi_{\eps,N} (q_1) =  
\Gamma ( q_1 )  + {\mathcal R}_{\eps,N}(q_1) 
\equiv \Big( {\mathcal A}(q_1) - 
\I a_{2d-1}(\fhi_1) \frac{q_1^{2d}}{2d}  \Big) 
+ {\mathcal R}_{\eps,N}(q_1) 
\end{equation}
where the quadratic form 
$$
{\mathcal A} (q) :=  \I \frac{( 2 + \eps )}{2}( \partial_{\fhi_2} q)^2 + 
( \partial_{\fhi_1} q )(\partial_{\fhi_2} q ) 
$$
is positive definite on $ Q_+ $, 
negative definite on $ Q_- $ and  
zero-definite on $Q_0 $. For $q_1 = q_+ + q_0 + q_- \in Q_1 $, 
\begin{equation}\label{qua}
{\mathcal A}(q_1) = {\mathcal A}(q_+ + q_0 + q_- ) = 
{\mathcal A}(q_+) + {\mathcal A}(q_- ) =
\frac{\alpha_+}{2} |q_+ |_{H^1}^2 - \frac{\alpha_{-}}{2} |q_- |_{H^1}^2 
\end{equation}
for suitable positive constants $ \alpha_+ $, $ \alpha_- $, bounded away
from $0$ by constants independent of $ \eps $.

We shall prove the existence of critical points 
of $ \Phi_{\eps,N} $ in $ B_{2R} $ 
 of ``linking type''.
\end{subsection}

\medskip 

\begin{subsection}{Linking critical points of the 
reduced action functional  $ \Phi_{\eps,N} $ }\label{lin}
We can not directly apply the linking 
Theorem because $ \Phi_{\eps,N} $ 
is defined only in the ball 
$ B_{2R} $. 
Therefore our first step is to 
extend $ \Phi_{\eps,N} $ to the whole space $ Q_1 $.
\\[1.5mm]
{\bf Step 1:} {\sl Extension of $ \Phi_{\eps,N} $.}
We define the extended action functional
$ {\widetilde \Phi}_{\eps,N} \in C^1 ( Q_1 , \mathbb{R}) $ as
$$
\Tf ( q_1 ) 
:=  \Gamma ( q_1 )  + {\widetilde{\mathcal R}}_{\eps, N } (q_1) 
$$
where $ \widetilde{\mathcal R}_{\eps, N}: Q_1 \to \R $ is  
$$
\widetilde{\mathcal R}_{\eps, N} (q_1) := 
\lam \Big( \frac{\h{q_1}^{2}}{R^{2}} \Big) {\mathcal R}_{\eps, N} (q_1)  
$$
and $ \lam : [0, +\infty ) \to [0,1] $ is a smooth, 
non-increasing, cut-off function such that 
$$ 
\left\{\begin{aligned}\lam(x) = 1 \qquad & |x|\leq 1 \\ 
\lam(x) = 0 \qquad & |x|\geq 4 \end{aligned}\right. \qquad |\lam'(x)| < 1 \, . 
$$ 
By definition $ {\widetilde \Phi}_{\eps,N} \equiv \Phi_{\eps,N}  $
on $ B_R := \{  q_1 \in Q_1 \, : \, |q_1|_{H^1} \leq R  \}$ and 
$ {\widetilde \Phi}_{\eps,N} \equiv \Gamma $ outside $ B_{2R} $.  

Moreover, by \refa{Restun}-\refa{Rest1un}, 
there is a constant $C_3 (R) \geq C_2 (R) > 0 $ such that
$ \forall |q_1|_{H^1} \leq 2 R $  
\begin{eqnarray}\label{Rest}
| {\widetilde{\mathcal R}}_{\eps,N} (q_1) | & \leq & 
C_3(R) \Big(\delta+ |\eps|\gam^{-1} + \frac{1}{N^2} \Big) \\
\label{Rest1}
\Big| {\widetilde{\mathcal R}}_{\eps,N}' (q_1)[h] \Big|  & \leq &
C_3 (R) \Big(\delta+ |\eps|\gam^{-1} + \frac{1}{N^2} \Big) | h |_{H^1}, 
\qquad \forall h \in Q_1 \, . 
\end{eqnarray}

In the sequel we shall always assume
$$
C_3 ( R ) \Big( \delta+ |\eps|\gam^{-1} + \frac{1}{N^2} \Big) \leq 1 \, .
$$

\medskip

\np
{\bf Step 2:}{ \sl $\Tf $ verifies the geometrical hypotheses 
of the linking Theorem.}

\fig{center}{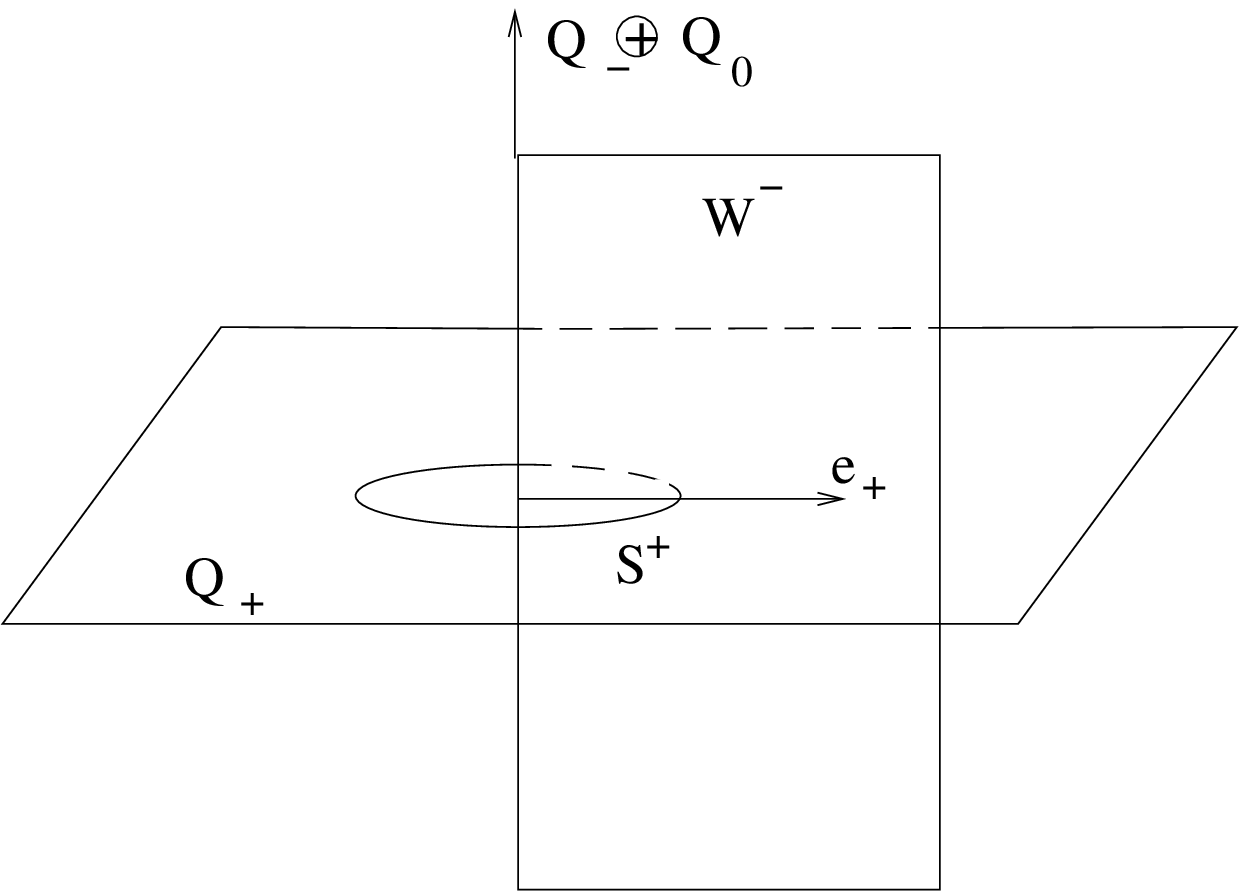}{The cylinder $W^-$ and the sphere $S^+$ link.}{0.7}
\medskip

\begin{lem}\label{L13}
There exist $ \eps $-$N$-$\gamma$-independent 
positive constants 
$ \rho $, $ \omega $, $ r_1 $, 
$ r_2 > \rho $, and $ 0 < $ $ \eps_1 (R) \leq $ $\eps_0 (R)  $, 
$ N_1 (R) \geq N_0 (R) $ 
such that, $ \forall |\eps| \gam^{-1} 
\leq \eps_1 (R)  $, $\forall N \geq N_1 (R) $
\\[0,6mm]

{\bf (i)} $\Tf(q_1)\geq  \omega > 0 $, $ \forall q_1 \in S^+ := 
\Big\{ q_1 \in Q_1 \cap Q_+ : \h{q_1}=\rho \Big\} $,

{\bf (ii)} 
$ \Tf ( q_1 ) \leq \omega \slash 2 $, $ \forall q_1 \in \partial W^- $ 
where $ W^- $ is the cylinder 
$$ 
W^-:= \Big\{ q_1 = q_0 + q_-+ r e^+ \,, \ \h{q_0 + q_-}\leq r_1\,, \  
q_- \in  Q_1 \cap Q_-,\; q_0\in \R \, , \, 
r \in [0,r_2]  \Big\} 
$$ 
and $ e_+ := \cos (\fhi_2 ) \in Q_1 \cap Q_+ $. Note that $\rho,\ome$ are {\it independent} of $R$.
\end{lem}
\np
In the following $ { \kappa }_i $, $ \kappa_\pm $
will denote positive constants 
{\it independent} on $ R $, $ N $, $ \eps $ and $ \gam $.
\begin{proof} {\bf (i)}  
$ \forall q_+ \in Q_1 \cap Q_+  $ with $ | q_+|_{H^1} = \rho < R $ we have
\begin{eqnarray}
{\widetilde \Phi}_{\eps,N} (q_+) = \Phi_{\eps,N} (q_+)& = & 
{\mathcal A} (q_+) - \I a_{2d-1} ( \fhi_1 ) \frac{q_+^{2d}}{2d} + 
{\mathcal R}_{\eps,N} (q_+) \nonumber \\ 
& \geq & \frac{\alpha_+}{2} \rho^2 - \kappa_1 \rho^{2d} 
- \Big( \delta+ |\eps|\gam^{-1} + \frac{1}{N^2} \Big) C_3(  R ) 
\label{dalb}.
\end{eqnarray}
Now we fix $ \rho > 0 $ small such that 
$ (\alpha_+ \rho^2 \slash 2) - 
\kappa_1 \rho^{2d} \geq \alpha_+ \rho^2 \slash 4 $. Next, for
$(\delta+|\eps|\gam^{-1} +$ $ N^{-2}) C_3(R) \leq$ 
$ \alpha_+ \rho^2 \slash 8 $  we get by (\ref{dalb})
$$
\Tf( q_+ ) \geq  \frac{\alpha_+}{8} \rho^2 =: \ome > 0, \qquad  
\forall \ q_+ \in Q_1 \cap Q^+ \quad {\rm with} \quad |q_+| = \rho \ .
$$

\np {\bf (ii)}  
Let 
\begin{eqnarray*}
B_1 & := & \Big\{ q_1 = q_0 + q_- + r_2 e_+ \ {\rm with} \ \ 
|q_0 + q_- |_{H^1} \leq r_1, q_- \in Q_1 \cap Q_- \Big\} 
\subset \partial W^- \\
B_2 & := & \Big\{ q_1 = q_0 + q_- + r e_+ \ {\rm with} \ \ 
|q_0 + q_- |_{H^1} = r_1, q_- \in Q_1 \cap Q_- ,
\ r \in [0, r_2 ]  \Big\} \subset \partial W^- 
\end{eqnarray*}
and choose $ r_1, r_2 > 2R $. For $ q_1 = q_0 + q_- + r e_+ \in B_1 \cup B_2 $ 
\begin{eqnarray}
{\widetilde \Phi}_{\eps,N} (q_1) = \Gamma ( q_1 ) & = & {\mathcal A}(q_1) 
- \I a_{2d-1}(\fhi_1) ( q_0 + q_- + r e_+ )^{2d} 
\nonumber \\
& = & - \frac{\alpha_-}{2} |q_-|_{H^1}^2 + r^2 {\mathcal A}(e_+) 
- \I a_{2d-1}(\fhi_1) \frac{( q_0 + q_- + r e_+ )^{2d}}{2d} \nonumber \\
& \leq & - \frac{\alpha_-}{2} |q_-|_{H^1}^2 + r^2 {\mathcal A}(e_+) 
- \alpha \I ( q_0 + q_- + r e_+ )^{2d} \label{ult10}
\end{eqnarray}
because $ a_{2d-1} (\fhi_1 ) \slash 2d  \geq \alpha > 0 $. 
Now, by H\"older inequality and orthogonality
\begin{eqnarray*}
\I ( q_0 + q_- + r e_+ )^{2d} & \geq & \kappa_2 
\Big( \I ( q_0 + q_- + r e_+ )^2 \Big)^{d}  \\
& = & \kappa_2 \Big( \I q_0^2 + q_-^2 + r^2 e_+^2 \Big)^{d} \\
& \geq & \kappa_3 ( q_0^2 + r^2 )^d \geq  \kappa_3 ( q_0^{2d} + r^{2d} ) 
\end{eqnarray*}
and by (\ref{ult10}) we deduce
$$
{\widetilde \Phi}_{\eps,N} ( q_0 + q_- + r e_+ ) \leq  
( \kappa_+ r^2 - \kappa_3 r^{2d} ) - \Big( \frac{\alpha_-}{2} |q_-|_{H^1}^2 +
\kappa_3 q_0^{2d} \Big) \ .
$$
Now we fix $ r_2 $ large such that $ \kappa_+ r_2^2 - 
\kappa_3 r_2^{2d} \leq 0 $
and therefore 
$$
{\widetilde \Phi}_{\eps,N} (q_1) 
\leq \kappa_+ r_2^2 - \kappa_3 r_2^{2d} \leq 0  \qquad \forall q_1 \in B_1 \ . 
$$
Next, setting $ M := \max_{ r \in [0,r_2]}(\kappa_+ r^2 - \kappa_3 r^{2d}) $,
we fix $ r_1 $ large such that 
$$
\frac{\alpha_-}{2} |q_-|_{H^1}^2 +
\kappa_3 q_0^{2d}  \geq M \qquad \forall \ |q_- + q_0 | = r_1 
$$
and therefore 
$$
{\widetilde \Phi}_{\eps,N} (q_1) 
\leq M - \Big( \frac{\alpha_-}{2} |q_-|_{H^1}^2 +
\kappa_3q_0^{2d} \Big) \leq 0 \qquad \forall q_1 \in B_2 \ . 
$$
Finally  if $ q_1 = q_- + q_0 $:
\begin{eqnarray}
\Tf (q_1)& = &  {\mathcal A}(q_- ) - 
\I a_{2d-1}(\fhi_1) \frac{q_1^{2d}}{2d} 
+ {\widetilde {\mathcal R}}_{\eps,N} (q_1) \nonumber \\
&\leq & |{\widetilde {\mathcal R}}_{\eps,N} (q_1)| 
\leq C_3(R )(\delta+ |\eps|\gam^{-1}
+  N^{-2}) \label{neg}
\end{eqnarray}
and so  $ \Tf (q_1) \leq \omega \slash 2 \ $ if 
$ \ C_3(R )(\delta+ |\eps|\gam^{-1} +  N^{-2}) \leq \omega \slash 2 $.
\end{proof}

\medskip

We introduce the minimax class 
$$ 
{\mathcal S} := \Big\{ \psi \in C( \overline{W^-}, Q) \ \  |  \ \   
\psi = {\rm Id}  \ \ {\rm on} \ \ \partial W^- \Big\} \, . 
$$
\np
The maps of $ {\mathcal S} $ have an important intersection property,
see e.g. Proposition 5.9 of \cite{ra}.
\begin{pro}\label{linkf}
{\bf($ S^+ $ and $ W^- $ link with respect to $ {\mathcal S} $.)} 
$$
\psi \in {\mathcal S} \quad \Longrightarrow \quad
\psi (W^-) \cap S^+ \neq \emptyset \, .
$$
\end{pro}

Define the minimax linking level 
$$ 
{\mathcal K}_{\eps,N} := 
\inf_{\psi\in {\mathcal S}}\max_{q \in W^-} \Tf (\psi(q_1))  \, .
$$
By the intersection property of Proposition \ref{linkf}
and Lemma \ref{L13}-($i$) 
$$
\max_{q_1 \in W^-} \Tf (\psi(q_1)) \geq 
\min_{q_1 \in S^+} \Tf ( q_1 ) \geq \omega > 0 \qquad 
\forall \psi \in {\mathcal S}  
$$
and therefore 
$$
{ \mathcal K }_{\eps,N} > \omega > 0 \, .
$$
Moreover, since Id $\in {\mathcal S} $ and (\ref{Rest})
\begin{eqnarray}
{\mathcal K}_{\eps,N} &\leq & \max_{q_1 \in W^-} \Tf (q_1)  \leq 
\max_{q_1 \in W^-} \Big( \Gamma ( q_1 ) 
+ {\widetilde {\mathcal R}}_{\eps,N} (q_1)\Big)  
\nonumber 
\\
& \leq & \max_{q_1 \in W^-} \Big( \frac{\alpha_+}{2} |q_+ |_{H^1}^2 + 
\frac{\alpha_{-}}{2} |q_- |_{H^1}^2  
+ \I \kappa q_1^{2d}\Big) + 1 
\leq {\mathcal K}_{\infty} < + \infty
\label{keps}
\end{eqnarray}  
where  ${\mathcal K}_{\infty}$ is independent
of $ N ,\eps, \gamma $ since the constants $r_1, r_2 $ in the definition 
of $ W^{-}$ are independent
of $ N ,\eps, \gamma $. 
 
We deduce, by the linking theorem 
the existence of a (Palais-Smale) sequence $ (q_j) \in Q_1 $
at the level $ { \mathcal K }_{\eps, N } $, namely 
\Beq{PS}{
\Tf (q_j)\rightarrow {\mathcal K}_{\eps,N} \, , 
\qquad \Tf' (q_j) \rightarrow 0 \ .}

\medskip 

\np
{\bf Step 3:}{ \sl Existence of a nontrivial critical point.}
Our final aim is to prove that 
the Palais-Smale sequence $ q_j $
converges, up to subsequence, to some non-trivial critical point 
$ {\overline q}_1 \neq 0 $ 
in some open ball of $ Q_1 $ where 
$ \Tf $ and $ \Phi_{\eps,N} $ coincide.

\begin{lem} \label{L28}
There exists a constant  $ R_* > 0 $, 
independent on $R$-$\eps$-$N$-$\gamma$, and functions 
$ 0 < \eps_2 (R) \leq \eps_1 (R) $, $ N_2 (R) \geq N_1 (R) $  
such that for all $ |\eps| \gamma^{-1} \leq \eps_2 (R) $, $ N \geq N_2 (R) $  
the functional $ \Tf  $ possesses a non-trivial  
critical point $ \overline{q}_1 \in Q_1 $
with critical value $ \Tf ( \overline{q}_1   ) = 
{\mathcal K}_{\eps,N}$, satisfying 
$ \h{ \overline{q}_1} \leq R_* $.
\end{lem}
\begin{proof}
Writing $\Tf (q) = \Gamma (q) + 
{\widetilde {\mathcal R}}_{\eps,N} (q) $ we derive,
by (\ref{Rest})-(\ref{Rest1})
\begin{eqnarray*}
\Tf ( q_j ) - \frac{1}{2} \Tf'(q_j ) [ q_j ] & = & 
\Gamma (q_j) -  \frac{1}{2} \Gamma'( q_j )[ q_j ] + 
\Big( {\widetilde {\mathcal R}}_{\eps,N} (q_j) - 
\frac{1}{2} {\widetilde {\mathcal R}}_{\eps,N}' ( q_j )[ q_j ] \Big)  \\
& = & \Big( \frac12 - \frac{1}{2d} \Big)\I a_{2d-1}(\fhi_1 ) q_j^{2d} +
\Big( {\widetilde {\mathcal R}}_{\eps,N} ( q_j ) 
- \frac{1}{2} {\widetilde {\mathcal R}}_{\eps,N}' ( q_j )[ q_j ] \Big) \\
& \geq & \alpha \Big( \frac12 - \frac{1}{2d} \Big)   \I q_j^{2d} -
(\delta+ |\eps|\gam^{-1} + N^{-2})  C_3 ( R ) \ .
\end{eqnarray*}
Therefore, by (\ref{keps})-(\ref{PS}) 
\Beq{bps}{ {\mathcal K}_\infty + 1 + \h{q_j} \geq \kappa_1 \I q_j^{2d} :=
\kappa_1 |q_j|_{L^{2d}}^{2d} \, .}
We also deduce, by (\ref{bps}), H\"older inequality and
orthogonality
 \begin{eqnarray}
 {\mathcal K}_\infty + 1 + \h{q_j} 
& \geq & \kappa_2 \Big( \I 
\Big( q_{+,j} + q_{0,j} + q_{-,j} \Big)^2 \Big)^d \nonumber \\
& = & \kappa_2 \Big( \I  
q_{+,j}^2 + q_{0,j}^2 + q_{-,j}^2 \Big)^d 
\geq \kappa_3  \ (q_{0,j})^{2d}  \nonumber 
\end{eqnarray}
and therefore 
\Beq{q0n}{
| q_{0,j} | \leq \kappa_4 \Big(  1 + \h{q_j} \Big)^{1\slash 2d} \ .}

By (\ref{Rest})-\refa{Rest1} and H\"older inequality
\begin{eqnarray}
\Tf'(q_j ) [ q_{+,j} ]  & = & 
\alpha_+ \h{q_{+,j}}^2 - \I a_{2d-1}( \fhi_1 ) q_j^{2d-1} q_{+,j} 
+ {\widetilde {\mathcal R}}_{\eps,N}' ( q_j )[ q_{+,j} ] \nonumber \\
& \geq  & \alpha_+ \h{q_{+,j}}^2 - \kappa_5 \h{q_{+,j}} \I |q_j|^{2d-1} -
(\delta+\gam^{-1}|\eps|+ N^{-2}) C_3(R) \h{q_{+,j}} \nonumber \\
& \geq  & \kappa_6 
\h{q_{+,j}} \Big( \h{q_{+,j}} -  |q_j|_{L^{2d}}^{2d-1} - 1 \Big) \ .
\label{poie}
\end{eqnarray}
By (\ref{poie}) and (\ref{bps}), using that $\Tf'(q_j ) \to 0 $ and 
simple inequalities, we conclude 
$$
\h{q_{+,j}} \leq \kappa_7 \Big( 1+ 
\h{q_j}^{(2d-1)\slash {2d}} \Big) \ .
$$
Estimating analogously $ \Tf'(q_j ) [ q_{-,j} ] $ we derive
$$
\h{q_{-,j}} \leq \kappa_8 \Big( 1+ \h{q_j}^{(2d-1)\slash {2d}} \Big) 
$$
and by (\ref{q0n}) we finally  deduce 
$$
\h{q_j} = |q_{0,j}| + \h{q_{+,j}} + \h{q_{-,j}} 
\leq \kappa_9 \Big( 1+  \h{q_j}^{1\slash {2d}} + 
\h{q_j}^{(2d-1)\slash {2d}} \Big) \ . 
$$
We conclude that
$ \h{q_j} \leq R_* $ for a suitable positive constant $ R_* $ independent 
of $ \eps $, $ N $, $ R $ and $ \gamma $. 

Since $ Q_1 $ is finite dimensional
$ q_j $ converges, up to 
subsequence, to some critical point $ \overline{q}_1 $ 
of $\Tf $ with $\h{\overline{q}_1} \leq R_*$. Finally, 
since $ \Tf (\overline{q}_1) =$ 
${\mathcal K}_{\eps, N} \geq $ $\ome > 0 $ we conclude that
$ \overline{q}_1 \neq 0 $. 
\end{proof}

We are now ready to prove Theorem A in the case $ \ome_1 = 1 $.

\begin{proof}[Proof of Theorem A for $ \ome_1 = 1 $]

Let us fix
$$
\bar R := R_* + 1 \qquad {\rm and \  take } \qquad 
|\eps| \gam^{-1} \leq \eps_2 ( \bar R ):= \overline{\eps} \, .  
$$
Set $\bar {N} := N_2 ( \bar{R} ) \geq N_0 (\bar{R})$.

Applying Lemma 2.3 we obtain, for
$$
0 < \sig \leq \frac{1}{N_2(\bar {R})} 
$$
a solution  
$ (q_2(q_1),p(q_1)) \in $ $ (Q_2 ( \bar{N} ) \oplus P) \cap
\B_{\sig,s} $ of the 
$(Q_2)$-$(P)$ equations $ \forall \h{q_1} \leq 2 \bar {R} $.
By Lemma \ref{L28} the extended 
functional $\widetilde\Phi_{\eps,N}(q_1)$ possesses
a critical point $ \bar q_1 \neq 0 $ with 
$ \h{\bar q_1} \leq R_*< \bar {R}  $. 
Since $ \widetilde \Phi_{\eps,N}(q_1)$ coincides with 
$\Phi_{\eps,N}(q_1) $ on the ball $ B_{\bar {R}} $
we get, by Lemma \ref{L24},   
the existence of a nontrivial weak solution 
$ \overline{q}_1 +q_2( \overline{q}_1) + p(\overline{q}_1) 
\in \B_{\sig,s} $ of equation (\ref{l1}). 
Finally 
$$
u = |\eps |^{1\slash 2(d-1)} \Big[  
\overline{q}_1 +q_2( \overline{q}_1) + p(\overline{q}_1) \Big]
\equiv |\eps |^{1\slash 2(d-1)} \Big[  
\overline{q}_\eps  + p(\overline{q}_1) \Big]
$$
solves equation \refa{e01}. 

The solution $ {\overline q}_\eps  := \overline{q}_1 + q_2( \overline{q}_1) $
of the $(Q)$-equation belongs to $ Q \cap \B_{\sig,s+2} $
by the regularizing properties
of $L_1^{-1} $, see in Lemma \ref{L1inv} formula 
\refa{L1-1}. 

Since $ \bar p := p(\overline{q}_1) $ solves
\begin{equation}\label{regu}
\Big( \partial^2_{\fhi_1} + 
2(1+\eps) \partial_{\fhi_1}\partial_{\fhi_2}\Big) \bar p   
= - \eps 
\Big[ (2 + \eps ) \partial^2_{\fhi_2} \bar p 
+ \Pi_P f (\fhi_1, u , \delta) \Big]
\in \B_{\sig',s} \quad \forall 0 < \sig' < \sig 
\end{equation}
and the eigenvalues of 
$ \partial^2_{\fhi_1} + 2(1+\eps) \partial_{\fhi_1}\partial_{\fhi_2} $
restricted to $ P $ 
satisfy, for $\eps \in {\mathcal B}_\gam $, 
$$
\Big| l_1 \Big[ (l_1 + 2 l_2) + \eps 2 l_2 \Big] \Big| 
\geq \gamma \frac{| l_1 |}{2|l_2 |} \qquad \forall \,
l_1 + 2 l_2 \neq 0, \, l_2 \neq 0  
$$
and 
we deduce that $ \bar p \in  \B_{\sig'',s+1}$ for all
$0 < \sig'' < \sig' $ and 
$ | \partial_{\fhi_1} \bar p |_{\sig'',s} = O( |\eps |\gam^{-1} )$. 
Now, again by \refa{regu}, 
$$
\partial^2_{\fhi_1} \bar p =  
-2(1+\eps) \partial_{\fhi_2}\partial_{\fhi_1} \bar p  
- \eps \Big[(2 + \eps ) \partial^2_{\fhi_2} \bar p 
+ \Pi_P f (\fhi_1, u , \delta)\Big] 
\in \B_{\overline{\sig},s} \quad \forall 0 < \overline{\sig} < \sig''\, . 
$$
therefore $ \bar p
\in \B_{\overline{\sig}, s+2} $ and 
$ | \bar p |_{\overline{\sig}, s+2} = O(|\eps| \gamma^{-1} )$.
\refa{ba12} follows with $\bar s := s+2 > 2 $. 
\\[2mm] 
\indent
By \refa{e00}, the function $ v(\eps, t, x )=$ $  
u (\eps, t, x+(1+\eps)t) $ is a 
solution of equation \refa{e-0} with $\ome_1 = 1 $. 
The frequency 
$ \ome_2 = 1 + \eps \notin \Q $ since $\eps \in {\mathcal B}_\gam $. 
To show that $ v(\eps, t, x ) $  
is quasi-periodic it remains to prove that 
$ u $     
depends on both the variables 
$ (\fhi_1, \fhi_2) $ independently.

We claim  that  $ {\bar q}_1 \notin Q_0 \oplus Q_- $, i.e. 
$  \overline{q}_+ (\fhi_2 ) \in $ $ Q_+\setminus \{0\}$, and therefore  
$ u $ depends on $ \fhi_2 $.  
Indeed by Lemma 2.6 we know that 
$ \tilde\Phi_{\eps,N}(\bar q_1) > \ome > 0 $
and $ | \overline{q}_1 |_{H^1( \T )} < \overline{R} $. 
On the other hand, by \refa{neg} in Lemma 2.5   
$ \tilde \Phi_{\eps,N}(q_- + q_0) < \ome \slash  2 $, 
for all $ |q_- + q_0 |_{H^1} \leq  \bar{ R }$, 
so that necessarily $ \bar q_1 \notin Q_0 \oplus Q_- $. 

We claim that any solution $ u $  
of \refa{l1} depending only on $ \fhi_2 $, namely solving 
\begin{equation}\label{onev}
(2+\eps) u''(\fhi_2 ) +  f(\fhi_1, u(\fhi_2),\delta)= 0 \, , 
\end{equation}
is $u(\fhi_2) \equiv 0$ .
Indeed,  by definition,  
$$
\delta^{2(d-1)}f(\fhi_1, u,\delta)= 
f(\fhi_1,\delta u) =
\sum_{k=2d-1}^\infty a_k(\fhi_1)(\delta u)^k $$ (recall sign$(\eps) = 1$). 
Consider now a smooth zero mean function $g(\fhi_1)$ such that 
$\int_0^{2\pi} a_k(\fhi_1)g(\fhi_1)\neq 0$ for some $k$
(recall that by assumption ({\bf H}) some of the $ a_k (\fhi_2) $ 
are not constant). By \refa{onev} we have 
$$
(2+\eps) u''(\fhi_2 )\int_0^{2\pi} g(\fhi_1) 
d\fhi_1 +
\int_0^{2\pi} f(\fhi_1, u(\fhi_2),\delta) g(\fhi_1)  d\fhi_1 = 0
$$ 
which implies, by the assumption {\bf (H)} on $f$, 
\begin{equation}\label{sequ}
\sum _{k=2d-1}^\infty 
[\delta u(\fhi_2)]^k \int_0^{2\pi} a_k(\fhi_1)g(\fhi_1) 
d\fhi_1 = 0  \, .
\end{equation}
The function $ G( z ) :=  \sum _{k=2d-1}^\infty 
b_k z^k $ with $ b_k := \int_0^{2\pi} a_k(\fhi_1)g(\fhi_1) d  \fhi_1$   
is a nontrivial analytic function. Therefore equation \refa{sequ}, i.e.  
$ G(\delta u(\fhi_2 ) ) = 0 $,
cannot have a sequence of zeros accumulating to $ 0 $.
So, for $\delta $ small enough, $u (\fhi_2 ) \equiv 0 $.
\end{proof}

\medskip 

\begin{proof}[Proof of Theorem A for any rational frequency 
$ \ome_1 = \frac{n}{m} \in \Q$]
Consider now equation \refa{e1} with $ \ome_1 = n \slash m $ 
where $ n, m $ are coprime integers.

The space $ Q $, formed by the solutions of 
$ \partial_{\fhi_1}(\frac{n}{m}\partial_{\fhi_1}+2 \partial_{\fhi_2})q = 0 $
can be written as
$$
Q  =   \Big\{ q = \sum_{l \in \mathbb{Z}^2 }
\hat q_l e^{\ii l \cdot \fhi } \in \B_{\sig,s} \  | \ \ 
\hat q_l = 0 \ \ \ {\rm for } \ \ \ l_1 (nl_1 + 2m l_2) \neq 0 \Big\} 
$$
and is composed by functions of the form
$$
q(\fhi)=q_+(\fhi_2)+q_-(2m\fhi_1-n\fhi_2)+q_0 \, .  
$$
Let $ P $ be the supplementary space to $ Q $ 
and perform the Lyapunov-Schmidt decomposition
like in \refa{q1}-\refa{q2}-\refa{p}.

For $ \eps $ in the Cantor set $ \Ba_\gam $,  
the eigenvalues 
$$
D_l = \Big( \frac{n}{m} l_1 +\eps l_2  \Big)  \, \Big( 
\frac{n}{m} l_1 + 2 l_2 + \eps l_2 \Big)
$$
of the linear operator $ {\mathcal L}_\eps $ 
can be bounded, arguing as in Lemma \ref{L1}, by
$$
|D_l| = \frac{|(n l_1+\eps m l_2)(n l_1+2 m l_2 + \eps m l_2)|}{m^2} 
> \frac{\gam}{m^2} \qquad  
\forall l_1 \neq 0, \; n l_1 + 2m l_2 \neq 0 \, . 
$$
As a consequence 
$$
\quad \Big| {\mathcal L}_\eps^{-1}[h] \Big|_{\sig,s} \leq 
\frac{m^2|h|_{\sig,s}}{\gamma}, \qquad \forall h \in P \, , 
$$ 
and, in solving the $(Q_2)$-$(P)$ equations as in   
Lemma \ref{L2}, we obtain the new restriction 
$$
\gam^{-1} |\eps | \leq  \frac{\eps_0(R)}{m^2}\,,\qquad N\geq N_0(R)
$$ 
and the bound (compare with \refa{co1})  
$ | p(q_1) |_{\sig,s} \leq C_0 (R)|\eps| \gamma^{-1} m^2.   $

The corresponding reduced action functional 
has again the form \refa{r1}-\refa{redf} with 
the different quadratic part  
$$
{\mathcal A}(q_1) = {\mathcal A}(q_+ + q_0 + q_- ) = 
{\mathcal A}(q_+) + {\mathcal A}(q_- ) =
\frac{\alpha_+}{2} |q_+ |_{H^1}^2 - 
n^2 \frac{\alpha_{-}}{2} |q_- |_{H^1}^2 
$$  
and therefore it still possesses a linking critical point
$\bar q_1 \in Q_1 $.

To prove the bound \refa{ba12} note that the
eigenvalues of 
$ \ome_1^2 \partial^2_{\fhi_1} + 2 \ome_1 (1+\eps) 
\partial_{\fhi_1}\partial_{\fhi_2} $ ($ \ome_1 = n \slash m $)
restricted to $ P $ 
satisfy, for $\eps \in {\mathcal B}_\gam $, 
$$
\ome_1 \Big| l_1 \frac{(n l_1 + 2 l_2m ) + \eps 2 l_2 m}{m} \Big| 
\geq  \frac{\ome_1 | l_1 |\gamma}{2 |l_2 |m^2} \qquad \forall \, 
nl_1 + 2m l_2 \neq 0, \, l_2 \neq 0  
$$
and therefore 
$ \bar p \in \B_{\overline{\sig}, s+2} $ and 
$ | \bar p |_{\overline{\sig}, s+2} = O(| \eps | m^2 
\slash \ome_1^3 \gamma )$.
\end{proof}

\medskip

\end{subsection}
\end{section}\begin{section}{{\bf Case B: 
$ \ \ome_1 \not\in \mathbb{Q} \ $}} \label{caseB}

We now look for solutions 
of equation (\ref{e1}) when the forcing frequency 
$ \ome_1 $  is an irrational number. 

To fix notations we shall prove Theorem B when 
$ \int_0^{2\pi} a_{2d-1}(\fhi_1)d\fhi_1>0$ and therefore $ \eps>0$, 
i.e. sign$(\eps) = 1$ . 

Fixed $ 0 < \overline{s} < s - 1 \slash 2 $, 
the Nemitskii operator
$ u \to f(\fhi_1, u , \delta ) 
\in $ $ C^\infty  ( B_\rho, \B_{\sig, \overline{s}}) $ 
since, if $ a_k ( \fhi_1 ) \in H^s (\T) $, then $ a_k ( \cdot ) 
\in {\mathcal H}_{\sig,\overline{s}}$, $ \forall \sig > 0 $,
$ 0 < \overline{s} < s -  1 \slash 2 $.

\smallskip

For $ \eps = 0 $ equation (\ref{e1}) reduces to 
\begin{equation}\label{0bif}
\ome_1 \partial_{\fhi_1} \Big( \ome_1 \partial_{\fhi_1} 
+ 2 \partial_{\fhi_2} \Big ) q = 0 
\end{equation}
and its solutions $ q $ 
form,  by the irrationality of 
$ \ome_1 $, the infinite dimensional subspace 
\begin{equation}
Q :=\Big\{ q  
\in \B_{\sig,\overline{s}} \ : \ \partial_{\fhi_1} 
q\equiv 0 \Big\} 
 =   \Big\{ q= q( \fhi_2 )\in H^{\sig}(\T) \Big\} \ . 
\end{equation}

To find solutions of 
(\ref{e1}) for $ \eps \neq 0 $, 
we perform a Lyapunov-Schmidt reduction
and we decompose the space 
$$
\B_{\sig,\overline{s}} = Q \oplus P 
$$ 
where $ Q \equiv H^{\sig}(\T) $ and 
\begin{eqnarray*}
P & := & 
\Big\{ p = \sum_{l \in \mathbb{Z}^2}
{\hat p}_l e^{\ii l \cdot  \fhi } \in \B_{\sig,\overline{s}}  \ | \ \  
{\hat p}_l = 0 \ \ \ {\rm for} \  \ \ l_1  = 0 \Big \}  \ .
\end{eqnarray*}
Projecting equation (\ref{e1}) onto the closed subspaces 
$ Q $ and $ P $, setting $u = q + p \in \B_{\sig,\overline{s}} $ with 
$ q \in Q $, $ p \in P $ we obtain  
\begin{eqnarray} 
(2+\eps) \ddot q +
\Pi_Q \Big[ f(\fhi_1,q + p, \delta ) \Big] & = & 0  \quad \quad\quad (Q) 
\label{Q} 
\\ 
{\mathcal L}_\eps [p] 
+ \eps \, \Pi_{P}  
\Big[ f(\fhi_1, q + p, \delta ) \Big] & = & 0 \quad \quad\quad 
(P) \label{P}
\end{eqnarray}
where $ \ddot q = \partial_{\fhi_2}^2 q  $, 
$ \Pi_Q : \B_{\sig,\overline{s}} \to Q $ is the  
projector onto $ Q $,    
$$
(\Pi_Q u)( \fhi_2 ) := \frac{1}{2\pi} 
\int_{0}^{2 \pi } u( \fhi_1, \fhi_2) \ d \fhi_1 
\, ,
$$
and $ \Pi_P = {\rm Id} - Q $ is the 
projector onto $ P $.  
\\[1mm]
\indent
We could proceed now as in the previous section performing
a finite dimensional reduction and applying variational methods. 
However, in this case, the infinite dimensional 
$(Q)$-equation 
can be directly solved by the Implicit Function Theorem
in a space of analytic functions. 

For this, it is 
useful to consider the  parameter $\delta $
(and $ \eps = \delta^{2(d-1)} $) in the right hand side of 
 (\ref{P}),
as an independent parameter $ \delta = \eta $, 
$ \eps  = \eta^{2(d-1)} $,   and to solve the equation   
\Beq{Pe}{
{\mathcal L}_\eps [p] 
+ \eta^{2(d-1)} \, \Pi_{P}  
\Big[ f(\fhi_1, q + p, \eta ) \Big] = 0 \qquad \quad\quad (P_\eta)}
for $ (\eps , \ome_1) $ in the Cantor set ${\mathcal C}_\gam $ 
and {\it for all} $ \eta $ small.
In this way we highlight the smoothness of the solution 
$ p (\eta , \eps, \cdot )$ of the ($P_\eta $)-equation 
(\ref{Pe}) in the variable $ \eta $.

\medskip

\begin{subsection}{Solution of the $(P_\eta)$-equation}
We first prove that the operator 
$ {{\mathcal L}_\eps} : P \to P $ 
has a bounded inverse 
when $ (\eps, \ome_1 ) $ belongs to the Cantor set 
$ \Ca_\gam $ 
defined in \refa{cgam}.


\begin{lem} \label{Cgam}
For any $ \eps_0 > 0  $ the Cantor set $ \Ca_\gam $ is uncountable.
\end{lem}

\begin{proof} 
Consider the set $\ov \Ca$ of couples $ x_1$, $ x_2\in \Ba_\gam$ such  that:
$$ 
x_1\in (-\eps_1,\eps_1 ) \,, \quad x_2\in \Big( 1+\eps_1,2-\eps_1 \Big) 
\,,\quad  x_1+x_2 \notin \Q\,,\quad x_1-x_2\notin \Q \ . 
$$ 
where $\eps_1=\eps_0\slash 2 $.
$\ov \Ca $ is an uncountable subset of $ \R^2 $ 
since for all $ x_1 \in \Ba_\gam $ the conditions $ x_1 \pm $ 
$ x_2 \notin \Q $ exclude only a countable set of values $ x_2 $. 
The Lemma follows since $ \Ca_\gam $ contains $ \psi^{-1} \ov \Ca $ where 
$ \psi : $ $ ( \eps, \ome_1) \to $ $ (\eps \slash \ome_1,$
$ (2+\eps) \slash \ome_1 )$ 
is an invertible map for $ (\eps,\ome_1) \in (-\eps_0,\eps_0)\times(1,2) $.
\end{proof}
The operator $ {\mathcal L}_\eps $ 
has eigenvalues
$D_l = $ $(\ome_1 l_1 +\eps l_2)(\ome_1 l_1 + 2l_2 + \eps l_2)$.

\begin{lem}\label{le0} 
For $ (\eps,\ome_1)\in \Ca_\gam $  
the eigenvalues $ D_l $ of $ {\mathcal L}_\eps $ restricted to $ P $ 
satisfy
\begin{equation}\label{dnb}
|D_l| = \Big| (\ome_1 l_1 +\eps l_2)(\ome_1 l_1 + 2l_2 + \eps l_2)
\Big| > \gam \ ,  \qquad  \forall  l_1 \neq 0 \, . 
\end{equation}
As a consequence, 
the operator $ {\mathcal L}_\eps : P \to P $ 
has a bounded inverse $ {\mathcal L}_\eps^{-1} $ 
satisfying 
\begin{equation}\label{lome} 
\Big| {\mathcal L}_\eps^{-1}[p] \Big|_{\sig,\overline{s}} 
\leq \frac{|p|_{\sig,\overline{s}}}{\gam }, \qquad \forall  p \in P \ .
\end{equation}
\end{lem} 
\begin{proof} 

Estimate (\ref{dnb}) is trivially satisfied if 
$ -l_1 \neq $ $ \frac{\eps}{\ome_1}  l_2 $ 
and $ -l_1 \neq $ $ \frac{2+\eps}{\ome_1} l_2 $. 
Now, if $ -l_1= [  \frac{\eps}{\ome_1}  l_2] $, then
$|(2+\eps)l_2 +$ $  \ome_1 l_1 | >$ 
$ | ( 2 + \eps ) l_2 - \eps l_2 | - \1 >$ $ |l_2|.$
Therefore, using $ |\ome_1 l_1 +\eps l_2| > \gam \slash | l_2 |  $, 
we get  (\ref{dnb}). 
The same estimate (\ref{dnb}) 
holds if  $ -l_1 =[  \frac{2+\eps}{\ome_1} l_2]$ since, 
in this case, $ |\ome_1 l_1 +\eps l_2| > $ 
$ |(2+\eps)l_2-\eps l_2|-\1> $ $ |l_2| $.
\end{proof}

\medskip

Fixed points of the nonlinear operator ${\mathcal G} : 
P \to P $ defined by 
$$
{\mathcal G} (\eta, p):= 
-\eta^{2(d-1)} \, {\mathcal L}_\eps^{-1} \Pi_P  f( \fhi_1 , q + p,\eta ) 
$$
are solutions of the $(P_\eta)$-equation.

\begin{lem}\label{le7} 
Assume $ (\eps, \ome_1 ) \in  \Ca_\gam $. 
$ \forall R > 0 $ 
there exists $ \eta_0 (R) $, $ C_0 (R) > 0 $ such that 
$ \forall  |q |_{H^{\sig}(\T)} \leq R $, $ 0 < \eta  \gamma^{-c} 
\leq \eta_0(R) $, with $c=1\slash 2(d-1)$,
there exists a unique $ p( \eta, q) \in 
P \cap {\mathcal H}_{\sig,\overline{s}}$  
solving the $(P_\eta)$-equation (\ref{Pe}) and satisfying 
\begin{equation}\label{pB}
|p  ( \eta, q )|_{\sig,\overline{s}} \leq C_0 (R) 
\eta^{2(d-1)} \gamma^{-1} 
\end{equation}
and the equivariance property 
\Beq{inva}{
p (\eta, q_\theta )(\fhi_1, \fhi_2 ) = p (\eta, q) (\fhi_1, \fhi_2 - \theta ),
\qquad \forall  \theta \in \T }
where $ q_\theta (\fhi_1, \fhi_2 ) := q (\fhi_1, \fhi_2 -\theta ) $.
Moreover $ p (\cdot, \cdot) \in C^1 ( (0, \eta_0(R)) \times Q ; P) $. 

\end{lem}

\begin{proof} In the Appendix. \end{proof}
\end{subsection}

\smallskip

\begin{subsection}{The $(Q)$-equation}

Once the $(P_\eta)$-equation has been solved by 
$ p(\eta, q) \in P $ there remains
the infinite dimensional bifurcation equation 
\Beq{Q2}{ 
(2+\eps)\ddot q+\Pi_{Q}\Big[ f(\fhi_1,q+ p (\eta, q),\eta) \Big] = 0\, .} 
Recalling \refa{fd}, 
the $(Q)$-equation (\ref{Q2}) 
evaluated at $ \eta = 0 $ 
reduces to the ordinary differential equation
\Beq{eqq}{ (2+\eps) \ddot q+  \langle a_{2d-1} \rangle  \, q^{2d-1}=0}
where 
$ \langle a_{2d-1} \rangle  := (1 \slash 2\pi )\int_0^{2\pi}
a_{2d-1}(\fhi_1) \, d\fhi_1  $. 

Equation (\ref{eqq}) is a superlinear autonomous 
Hamiltonian system with one degree of freedom and can be 
studied by a direct phase-space analysis. 

\begin{lem}\label{pro1} 
There exists $ \overline{\sig} > 0 $ such that, 
equation (\ref{eqq}) 
possesses a $ 2\pi $-periodic, analytic solution 
$ \overline{q} ( \fhi_2 )\in H^{\overline{\sig}}( \T ) $.
Morevoer, $ \overline{q} (\fhi_2 ) $ is non-degenerate
up to time translations, i.e. the linearized equation 
on $ \overline{q} $
\Beq{eqql}{(2+\eps) \ddot h +    (2d-1) \langle a_{2d-1} \rangle
\overline{q}^{2(d-1)} h =0}
possesses a one-dimensional space of $ 2 \pi $-periodic solutions,
spanned by $ \dot{\overline{q}} $.
\end{lem} 

\begin{proof}
Up to a rescaling, equation (\ref{eqq})
can be written as 
$ \ddot x = - V'(x) $ with potential energy  
$ V( x ) := x^{2d} $. All  solutions of such system  
are analytic and periodic with period 
$$ 
T(E)= 4 \int_{0}^{E^{\frac{1}{2d}}} 
\frac{dx}{\sqrt{2(E- x^{2d})}} = 
4 E^{\frac{1}{2d}- \frac12}\int_0^1 \frac{dx}{\sqrt{2(1-x^{2d})}} \, . 
$$
The equation $T(E)= 2 \pi $ has a solution 
$\bar q(\fhi_2)$ which is in 
$H^{\overline \sig}(\T)$ for some appropriate $\bar\sig > 0 $. 
The non-degeneracy of the corresponding 
$ 2 \pi $-periodic solution follows by
$$
\frac{dT}{dE}= 2 \Big( \frac{1}{d}- 1 \Big) 
E^{\frac{1}{2d}- \frac32}\int_0^1 \frac{dx}{\sqrt{2(1-x^{2d})}} \neq 0  
$$
and the next Proposition proved in the Appendix.
\begin{pro}\label{noniso}
Suppose the autonomous 
second order equation $ - \ddot x = V'(x) $, $ x \in \R $,  possesses a 
continuous family of non-constant 
periodic solutions $ x (E,t) $ with energy $ E $ and period $ T(E) $ 
satisfying the anysocronicity condition $ \frac{dT(E)}{dE} \neq 0 $. 
Then $ x (E,t) $ is non-degenerate up to time translations, i.e. 
the $T(E)$-periodic solutions of the linearized equation
\Beq{lin1}{- \ddot h = D^2 V ( x(E,t))h  }
form a one dimensional subspace spanned by $ (\partial_t x)(E,t) $.
\end{pro}
\end{proof}

From now on we fix $ \bar R 
:= |\overline{q}|_{H^{\overline{\sig}}(\T)} + 1 $ in Lemma \ref{le7}
and take $0 < \eta \gamma^{-c} \leq \eta_0 ( \bar R ) $.

By Lemma \ref{pro1} and \refa{inva}, 
we can construct solutions of the infinite dimensional 
bifurcation equation (\ref{Q2}) by means of the Implicit Function Theorem:

\begin{lem}\label{le8} There exist $ 0 < 
\eta_1 \leq \eta_0 (\bar R)$,
$ C_1 > 0 $ such that 
for all $ 0 < \eta\gam^{-c} \leq \eta_1 $,  
equation (\ref{Q2}) has a unique (up to translations) solution  
$ {\bar q}_\eta  (\fhi_2 ) \in H^{\overline{\sig}}(\T ) $ satisfying
$$
|{\bar q}_\eta  - \overline{q}|_{H^{\overline{\sig}}(\T)} \leq C_1 |\eta |\, .
$$
\end{lem}


\smallskip

\begin{proof}[Proof of Theorem B]
Setting again $ \delta \equiv \eta $,  
$ {\bar q}_\eps (\fhi_2 ) + p(\eps, {\bar q}_\eps  ) $ solves (\ref{e1}) and 
$$
u(\eps, \fhi ) = | \eps |^{1 \slash 2(d-1)} 
\Big[ {\bar q}_\eps (\fhi_2 ) + p(\eps, {\bar q}_\eps  ) \Big]  
$$
is a non trivial solution of \refa{e01}. The bound \refa{sB}
follows by \refa{pB}.
As in Theorem A the solution $ u $ 
depends on both the variables $(\fhi_1, \fhi_2 ) $.
Finally, the solution $ v(\eps, t , x ) := u(\eps, \ome_1 t ,x+\ome_2 t )$ of 
(\ref{e-0}) 
is quasi-periodic since, by the definition of $ \Ca_\gamma$, 
$ \ome_1 \slash \ome_2 =$ $ \ome_1 \slash (1+\eps) \notin \Q $. \end{proof}

\begin{oss}\label{othersol}
To prove existence of solutions of (\ref{e1}), i.e.
\refa{e-0}, 
it is sufficient that the second order equation (\ref{eqq})
possessess a continuous, nonisocronous 
family of non-constant periodic orbits
one of them having period $ 2 \pi \slash j $, see Proposition \ref{noniso}. 
\end{oss}

The hypothesys that the leading term in the nonlinearity $ f $ 
is an odd power of $ u $ is not of 
technical nature. The following non-existence result holds:

\begin{pro} {\bf (Non-existence)} \label{nonex} 
Let $ f(\fhi_1, u ) = a ( \fhi_1) u^D $ with 
$ D $ even and $ \int_0^{2\pi} a ( \fhi_1)$ $d \fhi_1 \neq 0 $.
$ \forall R > 0 $, there exists  $ \eps_0 > 0 $ such that  
$ \forall \sig \geq 0 $, $ \overline{s}> s- \1 $,
$ \forall (\eps, \ome_1) \in {\mathcal C}_\gamma $ with $|\eps | < \eps_0 $, 
equation (\ref{e01}) does not possess  
solutions $ u \in \B_{\sig,\overline{s}} $ in the ball 
$ | u |_{\sig,\overline{s}} \leq R | \eps |^{1\slash (D-1)} $.
\end{pro} 

\begin{proof}
We first rescale equation (\ref{e01}) with 
$ u \to | \eps |^{1\slash (D-1)} u $ obtaining 
\begin{equation}\label{risca}
{\mathcal L}_\eps u + |\eps | a(\fhi_1 ) u^D = 0 \, .
\end{equation}
Write any solution $ u_\eps \in B_{\sig,\overline{s}}(R) := \{  
u \in \B_{\sig,\overline{s}} \ : \ | u |_{\sig,\overline{s}} 
\leq R \} $ of \refa{risca} 
as $ u_\eps =$ $ q_\eps + p_\eps $ with $ q_\eps \in Q $, $ p_\eps \in P $.
$ p_\eps $ satisfies the $(P)$-equation 
$ {\mathcal L}_\eps p + |\eps | \Pi_P a(\fhi_1 ) u^D = 0 $ 
and therefore
$| p_\eps |_{\sig,\overline{s}} 
\leq C(R) |\eps | $. Then, for $\eps $ small enough,  
$ p_\eps \equiv  $ $ p (\eps, q_\eps ) $
where $ p(\eps, q_\eps ) $ is 
constructed as in Lemma \ref{le7} and 
satisfies the estimate $|p(\eps, q_\eps)|_{\sig,\overline{s}} \leq C |\eps | 
|q_\eps|_{H^{\sig} (\T)}^D $.

The projection $ q_\eps $ satisfies the $(Q)$-equation 
\Beq{lya}{
(2+\eps) \ddot q_\eps + 
{\rm sign}( \eps ) 
\Pi_Q \Big[ a(\fhi_1) (q_\eps + p ( \eps, q_\eps ) )^D \Big]  = 0 }
and therefore $|q_\eps |_{H^{\sig,2}(\T)} \leq C(R) $.

We claim that $q_\eps \to 0 $ in $H^{\sig}(\T) $ 
(and so in ${\mathcal H}_{\sig,\overline{s}}$)
for $\eps \to 0$. 
Indeed, from any subsequence $q_\eps $, we can extract 
by the compact embedding $ H^{\sig,2}( \T ) \hookrightarrow 
H^{\sig} ( \T )$ 
another convergent subsequence $q_{\eps_n}$ such that  
$ q_{\eps_n} \to \overline{q} \in H^{\sig} ( \T ) $. 
By (\ref{lya}), we deduce that
$$
2 \ddot{\overline{q}} + {\rm sign}(\eps) \langle a \rangle \overline{q}^D = 0
$$
where $ \langle a \rangle := \int_0^{2\pi} a (\fhi_1 ) \ d \fhi_1  \neq 0$. 
Such equation does not possess non-trivial periodic solutions for both 
sign$(\eps) = \pm 1 $, i.e. $ \eps > 0 $ and $ \eps < 0  $, and we conclude
that $\overline{q} = 0 $. 

We finally prove that equation (\ref{lya}) 
does not possess non-trivial periodic solutions 
in a small neighborhood of the origin.

Linearizing equation (\ref{lya}) at $ q = 0 $ we get 
$(2 +\eps )  \ddot h = 0 $
whose solutions in $ H^{\sig}(\T) $ are the constants. 
We can perform another Lyapunov-Schmidt reduction close to $ 0 $ 
decomposing $ H^{\sig}(\T) = $ $\{$constants$\}$ $ \oplus $ 
$\{$zero average functions$\}$, namely
$ q_\eps = \rho + w $.
By the Implicit function Theorem we get
that for any constant $ | \rho | \leq \rho_0 $ small enough
(independently of $\eps $) there exists a 
unique zero average function $ w_\rho $ with $|w_\rho |_{H^{\sig}( \T )}
= O(\rho^D )$ solving
$$
(2+\eps)\ddot w_\rho
+ \Big[ a(\fhi_1) ( \rho + w_\rho + p ( \eps, q_\eps ) )^D -
\Big\langle a(\fhi_1) (\rho 
+ w_\rho + p ( \eps, q_\eps ) )^D\Big\rangle \Big] = 0 \, .
$$
Hence $ \rho $ is such that 
$$
0 = 
\Big\langle a(\fhi_1) (\rho + w_\rho + p ( \eps, q_\eps ) )^D
\Big\rangle = 
\langle a \rangle \rho^D + o( \rho^D ) \, .
$$ 
This implies $ \rho = 0 $ 
since $ \langle a \rangle \neq 0 $ and so
$q_\eps = \rho + w_\rho = 0$. 
\end{proof}

\end{subsection}
\end{section}

\begin{section}{Appendix}

\begin{lem}\label{banach}
${\mathcal H}_{\sigma,s}$ is a Banach algebra for $ \sigma, s  \geq 0 $.
\end{lem}

\begin{proof}
By the product Cauchy formula 
$$
u v = \sum_{j \in \Z^2 } \Big( \sum_{k \in \Z^2} u_{j-k} v_k \Big)
e^{\ii j \cdot \fhi }
$$
and therefore
\begin{eqnarray*}
|u v |_{\sig,s} & := & 
\sum_{j \in \Z^2 } e^{\sig |j_2|} [j_1]^s 
\Big| \sum_{k \in Z^2} u_{j-k} v_k \Big| \leq 
\sum_{j \in \Z^2 } e^{\sig |j_2|} [j_1]^s 
\sum_{k \in \Z^2} |u_{j-k}| \, | v_k|  
\\
& \leq & \sum_{k \in \Z^2 } | v_k | 
\sum_{j \in \Z^2} |u_{j-k}| e^{\sig |j_2|} [j_1]^s \\
& \leq & 2^s \sum_{k \in \Z^2 } | v_k | e^{\sig |k_2|} [k_1]^s 
\sum_{j \in \Z^2} |u_{j-k}| e^{\sig |j_2-k_2|} [j_1 - k_1]^s 
:= 2^s |u|_{\sig,s} |v|_{\sig,s}
\end{eqnarray*}
since 
$e^{\sig |j_2|} \leq 
e^{\sig |j_2 - k_2|} e^{\sig |k_2|} $ and 
$[j_1] \leq 2 [j_1 - k_1] [k_1] $ for all $k, j \in \Z^2 $. 
\end{proof}

{\sl Proof of Lemma \ref{L2}.} 
Let us consider 
$$ 
B:=\Big\{ ( q_2 , p ) \in Q_2 \oplus P:  \ |q_2|_{\sig,s} 
\leq \rho_1\,, \ |p|_{\sig,s}\leq \rho_2 \Big\} 
$$
with norm $|( q_2 , p )|_{\sig,s} := |q_2 |_{\sig,s} + |p|_{\sig,s}$.
We claim that, under the assumptions \refa{hp}
there exists $ 0 < \rho_1 $, $ \rho_2 < 1 $, see \refa{rho12}, such that
the map $(q_2,p) \to \mathcal G(q_2,p; q_1) $ is a contraction in $B$, i.e.:
\begin{enumerate} \item[(i)] $(q_2,p)\in B \ \Longrightarrow $ 
$ \mathcal G(q_2, p; q_1 ) \in B$;
\item[(ii)] 
$
| \G(q_2,p; q_1)-\G(\widetilde q_2,\widetilde p; q_1) |_{\sig,s} \leq 
(1 \slash 2)  |(q_2,p) - (\widetilde q_2, \widetilde p) |_{\sig,s} \, , $
$ \ \forall (q_2,p), (\widetilde q_2,\widetilde p)\in B $.
\end{enumerate}
In the following ${\kappa}_i $ will denote positive constants 
{\it independent} on $ R $, $ N $ and $ \eps $ (i.e. on 
$ \delta := |\eps|^{1\slash 2(d-1)}$).

By \refa{l1invN} and the Banach algebra property of $ \B_{\sig,s} $
\begin{eqnarray} \label{co3}
| \G_1(q_2,p; q_1) |_{\sig,s} & = &   
| L_1^{-1}\Pi_{Q_2}  f(\fhi_1, q_1 + q_2 + p ,\delta) |_{\sig,s} \nonumber \\
& \leq & 
\frac{\kappa_1}{N^2} \Big( |q_1|_{\sig,s}^{2d-1}+ |q_2|_{\sig,s}^{2d-1} 
+ | p|_{\sig,s}^{2d-1} \Big )\,  
\end{eqnarray}
provided that $ 0 \leq \delta 
\leq \delta_0(R)$. 
Similarly, for $\eps \in {\mathcal B}_\gamma $, by \refa{b1},  
\begin{eqnarray} \label{co4}
| \G_2(q_2,p; q_1) |_{\sig,s} &= &   
| \eps {\mathcal L}_\eps^{-1} \Pi_P  f(\fhi_1 , q_1 + q_2 + p,\delta ) 
|_{\sig,s} \nonumber \\
& \leq &  
\kappa_2 | \eps | \gamma^{-1} \Big( |q_1|_{\sig,s}^{2d-1}+ |q_2|_{\sig,s}^{2d-1} 
+ | p|_{\sig,s}^{2d-1} \Big )\, . 
\end{eqnarray}
For all $ q_1 \in Q_1(N) $ and since $ 0 \leq s < 1\slash 2 $
\begin{eqnarray}
|q_1|_{\sig,s} & = &\sum_{|l_2| \leq N} 
|\hat q_{0,l_2}|e^{\sig |l_2|}  + 
| \hat q_{-2l_2,l_2}| e^{\sig |l_2|}[-2l_2]^s 
\nonumber \\
& \leq  & e^{\sig N} \sum_{|l_2| \leq N} 
|\hat q_{0,l_2}|  + | \hat q_{-2l_2,l_2}| [-2l_2]^s \nonumber 
\leq  \kappa_3 \Big[ \Big( \sum_{|l_2| \leq N} 
|\hat q_{0,l_2}|^2 [l_2]^2 \Big)^{1\slash 2}
\Big( \sum_{l_2 \in \Z }\frac{1}{[l_2]^2} \Big)^{1 \slash 2} \\
& + & \Big( \sum_{|l_2| \leq N} 
| \hat q_{-2l_2,l_2} |^2 [l_2]^2 \Big)^{1 \slash 2} 
\Big(\sum_{l_2 \in \Z } \frac{1}{ [l_2]^{2 (1-s)} } 
\Big)^{1 \slash 2} \Big] 
\leq \kappa_4 |q_1|_{H^1} \label{equiv}
\end{eqnarray}
whenever $ 0 \leq \sig N \leq 1 $. 

\np
Substituting in \refa{co3}-\refa{co4} we get
$ \forall | q_1 |_{H^1} \leq 2R $, $ \forall | q_2 |_{\sig,s} \leq \rho_1 
$, $ \forall | p |_{\sig,s} \leq \rho_2 $
\begin{eqnarray}
|\G_1(q_2, p ; q_1) |_{\sig,s} & \leq & 
\kappa_5 N^{-2} \Big( R^{2d-1} + \rho_1^{2d-1} 
+ \rho_2^{2d-1}  \Big) \label{r1a} \\
|\G_2(q_2,p; q_1 ) |_{\sig,s} & \leq & 
\kappa_5 | \eps | \gamma^{-1}\Big( R^{2d-1} + \rho_1^{2d-1} 
+ \rho_2^{2d-1}  \Big) \, \label{r2} .
\end{eqnarray}
Now, setting $ C_0 (R) :=  \kappa_5 R^{2d-1} $, we define 
\Beq{rho12}{
\rho_1 := \frac{2 C_0(R)}{N^2} \qquad \quad  \rho_2 := 2 |\eps|\gamma^{-1} 
C_0(R) \, .}
By \refa{r1a}, \refa{r2}
there exists $ N_0 (R) \in {\mathbb N}^+ $ and $ \eps_0(R) >0 $
such that $ \forall N \geq N_0 (R) $ and 
$ \forall | \eps | \gamma^{-1} \leq \eps_0 (R) $ 
$$
|\G_1( q_2, p ; q_1) |_{\sig,s} \leq \rho_1 
\qquad | \G_2 ( q_2, p ; q_1) |_{\sig,s} \leq \rho_2  
$$
proving (i).
Item (ii) is obtained with similar estimates.
 
By the Contraction Mapping Theorem  
there exists a unique fixed point $ (q_2 (q_1), p (q_1)) := $
$ (q_2 (\eps, N, q_1)$, $ p (\eps, N, q_1)) $
of $ \G $ in $ B $. 
The bounds \refa{co1} follow by \refa{rho12}. 

Since $ \G \in C^1( Q_2 \oplus P \times Q_1; 
Q_2 \oplus P \times Q_1 )$ 
the Implicit function Theorem implies that 
the maps $ Q_1 \ni 
q_1 \rightarrow (q_2 (\eps, N, q_1), p (\eps, N, q_1)) $ are $ C^1 $.   

Differentiating $ (q_2(q_1), p(q_1)) = 
\G (q_2(q_1), p(q_1), q_1) $ 
$$\begin{aligned} q_2'(q_1)[h] 
= - L_1^{-1}\Pi_{Q_2}  (\partial_u f)(\fhi_1, q_1 + q_2(q_1) + p(q_1),\delta )
\Big(h+q_2'(q_1)[h]+p'(q_1)[h] \Big) 
\\ p'(q_1)[h]= - \eps\La_\eps^{-1}\Pi_{Q_2}  
(\partial_u f)(\fhi_1, q_1 + q_2(q_1) + p(q_1),\delta )
\Big(h+q_2'(q_1)[h]+p'(q_1)[h] \Big)\end{aligned} 
$$  
and using \refa{l1invN}, \refa{b1} and the Banach algebra property
of $ \B_{\sig,s} $
$$ 
\begin{aligned}& |q_2'(q_1)[h] |_{\sig,s}\leq 
C(R)N^{-2} \Big( 
|h|_{\sig,s}+|q_2'(q_1)[h] |_{\sig,s}+|p'(q_1)[h] |_{\sig,s}\Big) 
\\ 
&|p'(q_1)[h] |_{\sig,s}\leq C(R)|\eps| \gamma^{-1} \Big( 
|h|_{\sig,s}+ | q_2'(q_1)[h] |_{\sig,s}+
| p'(q_1)[h] |_{\sig,s} \Big)
\end{aligned}
$$ 
which implies the bounds \refa{co2} since
$$
\det
\left|\begin{array}{lll} 
1-C(R)N^{-2}  & & -C(R) N^{-2} \\& & \\ -C(R)|\eps|\gamma^{-1} 
& & 1-C(R)|\eps| \gamma^{-1} \end{array}\right|\geq \1 
$$ 
for $ C(R)( |\eps|\gamma^{-1} +  N^{-2}) $ small enough and 
\refa{equiv}. \mbox

\bigskip

\np
{\sl Proof of Lemma \ref{L24}. }
By \refa{q2}, \refa{p} we have that, at 
$ u := q_1 + q_2(q_1) + p(q_1) $,
\begin{equation}\label{q2p}
d \Psi_\eps (u)[h] = 0 \, 
\ \forall h \in Q_2 \qquad {\rm and} \qquad
d \Psi_\eps (u)[h] = 0 \, \ \forall h \in P \, .
\end{equation}
Since $ q'_2 ( q_1 ) [k]\in Q_2 $ and 
$  p' (q_1)[k] \in P $ $ \forall k \in Q_1 $, we deduce
$$
d \Phi_{\eps,N} (q_1)[k] = d \Psi_\eps ( u)
\Big[ h + q'_2(q_1)[k] + p' (q_1)[k] \Big] = 
d \Psi_\eps ( u ) [k] \, \quad \forall k \in Q_1 
$$
and therefore $ u := q_1 + p(q_1)+q_2(q_1) $ solves also the 
($Q_1$)-equation (\ref{q1}). 

Write
$ \Psi_\eps (u) = \Psi^{(2)}_\eps(u) - \eps \I F(\fhi_1,u, \delta) $ 
where 
$$
\Psi^{(2)}_\eps(u):= \I \1 ( \partial_{\fhi_1} u )^2 + 
(1+\eps)( \partial_{\fhi_1} u) ( \partial_{\fhi_2} u ) + 
\frac{  \eps (2 + \eps) }{2} ( \partial_{\fhi_2} u )^2
$$ 
is an homogeneous functional of degree two.
By homogeneity:
\Beq{fu2}{\Psi_\eps (u)= \1 d\Psi^{(2)}_\eps (u)[u] 
- \eps \I F(\fhi_1,u,\delta)\, . }
By \refa{q2}, \refa{p} (i.e. \refa{q2p}) 
\begin{equation}\label{q2p2}
d\Psi^{(2)}_\eps (q_1+q_2(q_1) + p(q_1))[q_2(q_1) + p(q_1)]
= \eps \I f(\fhi_1,u, \delta)( q_2(q_1) + p(q_1) ) \, .
\end{equation}
Substituting in \refa{fu2} we obtain, at $u = q_1 + q_2(q_1) + p(q_1)$ 
\begin{eqnarray*} 
\Phi_{\eps,N}(q_1) & = & \Psi_\eps (q_1+p(q_1)+q_2(q_1))
=  \1 d\Psi^{(2)}_\eps (u)[ q_1+p(q_1)+q_2(q_1) ] 
- \eps \I F( \fhi_1 , u,\delta ) \\
& = & \1 d\Psi^{(2)}_\eps ( q_1 )[ q_1 ] 
- \eps \I F( \fhi_1 , u,\delta ) + 
\frac{1}{2} f(\fhi_1,u,\delta)( q_2(q_1) + p(q_1)) \\ 
& = & 
\Psi_0 (q_1) + \eps \I 
\frac{(2 + \eps) }{2} ( \partial_{\fhi_2} q_1 )^2
+  ( \partial_{\fhi_1} q_1) ( \partial_{\fhi_2} q_1 ) -
F(\fhi_1,u,\delta) \\
&+ & \1 f(\fhi_1,u,\delta) ( q_2 (q_1) + 
p(q_1)) = {\rm const} + \eps ( \Gamma (q_1) + {\mathcal R}_{\eps,N} (q_1))
\end{eqnarray*}
because $ \Psi_0 ( q_1 ) \equiv $ const. 

By \refa{co1} the bounds \refa{Restun}-\refa{Rest1un} follow.

\bigskip

\np
{\sl Proof of Lemma \ref{le7}. }
The existence of $ p(\eta, q) \in \B_{\sig,\overline{s}} $
can be proved as in Lemma 2.3
using the Contraction Mapping Theorem. The smoothness 
of $ p(\eta, q) $ follows by the Implicit Function Theorem
since $ {\mathcal G} (\eta, p)$  is smooth in $ \eta $ and $ q $.

By the invariance of equation (\ref{Pe}) under translations 
in the $ \fhi_2 $ variable 
the function $ p (\eta, q) (\fhi_1, \fhi_2 - \theta )$ solves
$$
p (\eta, q) (\fhi_1, \fhi_2 - \theta ) + 
\eta^{2(d-1)} \, {\mathcal L}_\eps^{-1} \Pi_P  f \Big( \fhi_1 , q_\theta + 
p (\eta, q) (\fhi_1, \fhi_2 - \theta ),\eta \Big) = 0 
$$
and, therefore, by uniqueness (\ref{inva}) holds.

\bigskip

\np
{\sl Proof of Proposition \ref{noniso}. }
Write $ x(E,t) =$ $ y( \ome(E) t , E ) $ 
where $ y ( \fhi , E ) $ is $ 2 \pi $-periodic in $ \fhi $ and 
$\ome (E) :=$ $2 \pi \slash T(E) $. The functions
$ (\partial_t x)(E,t) $ and 
\Beq{solli}{
(\partial_E x)(E,t) = t \frac{d\ome(E)}{dE} 
(\partial_\fhi y)( \ome(E) t , E ) + (\partial_E y) ( \ome (E)t , E)}
are two linearly 
independent solutions of the linearized equation (\ref{lin1}).  
$ (\partial_t x)(E,t) $ 
is $ 2\pi $-periodic while, since 
$$
\frac{d \ome(T)}{dT} = 2\pi T(E)^{-2} \frac{dE(T)}{dT}
\neq 0 \qquad {\rm and} \qquad (\partial_\fhi y) ( \fhi , E) \not \equiv 0
$$
(if not $ x(E,t) $ would be constant in $t$), $(\partial_E x)(E,t)$ is not 
$2\pi$-periodic. We conclude that the space of $T(E)$-periodic
solutions of (\ref{lin1}) form a $1$-dimensional linear space 
spanned by $ (\partial_t x)(E,t) $.

\end{section}

\end{document}